\documentclass[reqno,a4paper,oneside]{amsart} 
\usepackage{mathptmx,verbatim,cite,amsmath,url}

\newtheorem{theorem}{Theorem}[section]
\newtheorem{proposition}[theorem]{Proposition}
\newtheorem{lemma}[theorem]{Lemma}
\newtheorem{corollary}[theorem]{Corollary}
\newtheorem{example}[theorem]{Example}
\newtheorem{remark}[theorem]{Remark}
\newtheorem{definition}[theorem]{Definition}

\newcommand{\klg}[1]{\left\{#1\right\}}
\newcommand{\klr}[1]{\left(#1\right)}
\newcommand{\kle}[1]{\left[#1\right]}
\newcommand{\klrk}[1]{(#1)}
\newcommand{\klek}[1]{[#1]}
\newcommand{\klgk}[1]{\{#1\}}
\newcommand{\betrag}[1]{{\left|#1\right|}}

\newcommand{\norm}[1]{\left\|#1\right\|}

\def\sto{\to}
\def\supp{{\rm supp}}

\def\IR{{\mathbb{R}}}
\def\IN{{\mathbb{N}}}

\DeclareMathOperator{\Cap}{Cap}

\def\metric{\mathsf{M}}
\def\domain{\mathsf{D}}
\def\topo{\mathsf{T}}
\def\dist{{\mathsf{d}}}

\def\cD{\mathcal{D}}
\def\cF{\mathcal{F}}
\def\cL{\mathcal{L}}
\def\cM{\mathcal{M}}
\def\cN{\mathcal{N}}
\def\cP{\mathcal{P}}
\def\cT{\mathcal{T}}

\def\sC{\;\mathsf{C}}
\def\sF{\;\mathsf{F}}

\def\sW{\;\mathsf{W}}
\def\tsW{\mathcal{W}}

\def\tOm{\widetilde\Omega}

\def\su{\mathsf{u}}
\def\sg{\mathsf{g}}
\def\sv{\mathsf{v}}

\begin{document}

\title{Lattice Homomorphisms between Sobolev Spaces}
\author{Markus Biegert}
\dedicatory{\upshape%
  Institute of Applied Analysis, University of Ulm, 89069 Ulm, Germany\\
  \texttt{markus.biegert@uni-ulm.de}
}

\begin{abstract}
 We show in Theorem \ref{thm:zero} that every vector lattice homomorphism
 $T$ from $\sW^{1,p}_0(\Omega_1)$ into $\sW^{1,q}(\Omega_2)$ for
 $p,q\in (1,\infty)$ and open sets $\Omega_1,\Omega_2\subset\IR^N$ has a representation of
 the form $T\su=(\su\circ\xi)g$ ($\Cap_q$-quasi everywhere on $\Omega_2$)
 with mappings $\xi:\Omega_2\to\Omega_1$ and $g:\Omega_2\to[0,\infty)$.
 This representation follows as an application of an abstract and more
 general representation theorem (Theorem \ref{thm:art}).
 In Theorem \ref{thm:nonzero} we prove that every lattice homomorphism $T$
 from $\tsW^{1,p}(\Omega_1)$ into $\sW^{1,q}(\Omega_2)$
 admits a representation of the form
 $T\su=(\su\circ\xi)g$ ($\Cap_q$-quasi everywhere on $\Omega_2$)
 with mappings $\xi:\Omega_2\to\overline\Omega_1$ and $g:\Omega_2\to[0,\infty)$.
 Here $\tsW^{1,p}(\Omega_1)$ denotes the closure of $\sW^{1,p}(\Omega_1)\cap C_c(\overline\Omega_1)$
 in $\sW^{1,p}(\Omega_1)$ and every $\su\in\tsW^{1,p}(\Omega_1)$ admits a trace
 on the boundary $\partial\Omega_1$ of $\Omega_1$. Finally, in Theorem \ref{thm:nonzero2}
 we prove that every lattice homomorphism $T$ from $\tsW^{1,p}(\Omega_1)$ into $\tsW^{1,q}(\Omega_2)$
 where $\Omega_1$ is bounded has a representation of the form
 $T\su=(\su\circ\xi)g$ ($\Cap_{q,\Omega_2}$-quasi everywhere on $\overline\Omega_2$)
 with mappings $\xi:\overline\Omega_2\to\overline\Omega_1$ and $g:\overline\Omega_2\to[0,\infty)$. 
\end{abstract}

\maketitle

{\smaller[1]\tableofcontents}

%--------------------------------------------------------------------------------------------------
%--------------------------------------------------------------------------------------------------
%--------------------------------------------------------------------------------------------------
%--------------------------------------------------------------------------------------------------
%--------------------------------------------------------------------------------------------------

\section{Introduction}

Let $A_1$ and $A_2$ be non-empty sets, $E$ be a Banach space and let
$F_1(A_1,E)$ and $F_2(A_2,E)$ be linear spaces of $E$-valued functions
defined on $A_1$ and $A_2$, respectively. If $\xi:A_2\to A_1$ is such that
$u\circ\xi$ belongs to $F_2(A_2,E)$ for every $u\in F_1(A_1,E)$, then the operator $S_\xi$ which
maps $u$ to $u\circ\xi$ is called a composition transformation.
If $g:A_2\to\IR$ is a map such that $(u\circ\xi)g$ belongs to
$F_2(A_2,E)$ for every $u\in F_1(A_1,E)$, then the mapping
$T:u\mapsto (u\circ\xi)g$ is called a weighted composition
transformation induced by the CoMu-Representation $(\xi,g)$. There are
many results in the literature which assert that certain operators
are weighted composition transformations.

The classical {\em Banach-Stone Theorem} (Stefan Banach \cite{banach:32:tol}
and Marshall Harvey Stone \cite{stone:37:atb}) says the following.
Let $\topo_1$ and $\topo_2$ be
compact Hausdorff spaces. Then a bounded linear operator
$C(\topo_1)\to C(\topo_2)$ is a surjective isometry if and only if
$T$ has a CoMu-Representation $(\xi,g)$ for some
homeomorphism $\xi:\topo_2\to\topo_1$ and some continuous function
$g:\topo_2\to\klg{-1,1}$.
A simplified version of {\em Lamperti's Theorem} (John Lamperti \cite{lamperti:58:icf})
says the following. Let $1\leq p<\infty$, $p\not=2$ and let $T$ be a linear isometry of
$L^p([0,1])$ into itself. Then there is a Borel measurable mapping $\xi$
of $[0,1]$ onto (almost all of) $[0,1]$ and $\sg\in L^p([0,1])$ such that
$T\su=(\su\circ\xi)g$ for all $\su\in L^p([0,1])$.
Isometries on Orlicz spaces were considered by
John Lamperti in 1958 \cite{lamperti:58:icf} and by
G\"unter Lumer in 1963 \cite{lumer:63:orlicz}.
Interesting results for isometries between Sobolev spaces were obtained by
Geoff Diestel and Alexander Koldobsky in 2006 \cite{diestel:06:sti}
by considering $W^{1,p}(\Omega)$ as a subspace of a certain $L^p$ space.

In this article we consider vector lattice homomorphisms $T:\sW^{1,p}(\Omega_1)\to\sW^{1,q}(\Omega_2)$.
A large class of such lattice homomorphisms can be obtained as follows.
Let $\Omega_1$ and $\Omega_2$ be non-empty open sets in $\IR^N$. Then the
class $\cT^1_p(\Omega_1,\Omega_2)$ (defined in the book of
Maz'ya and Shaposhnikova \cite[Section 6.4.3]{mazya:85:tom})
consists by definition of those mappings $\xi:\Omega_2\to\Omega_1$ such that
$u\circ \xi\in \sW^{1,p}(\Omega_2)$ and
$\norm{u\circ\xi}_{\sW^{1,p}(\Omega_2)}\leq C\cdot \norm{u}_{\sW^{1,p}(\Omega_1)}$
for all $u\in\sW^{1,p}(\Omega_1)$, where $C$ is a constant independent of $u$.\
\footnote{For the case $p\in (N-1,N)$ see Gol'dshtejn and Romanov \cite[Section IV]{goldshtejn:84:tps}.}
A real-valued function $g$ defined on $\Omega_2$ belongs by definition to the class of Sobolev multipliers
(see Maz'ya and Shaposhnikova \cite[Chapter 1 and 6]{mazya:85:tom})
$\cM\klr{\sW^{1,p}(\Omega_2)\sto\sW^{1,q}(\Omega_2)}$ if $gu\in\sW^{1,q}(\Omega_2)$ for all
$u\in\sW^{1,p}(\Omega_1)$. Then $T:\sW^{1,p}(\Omega_1)\to\sW^{1,q}(\Omega_2)$ defined by
$Tu:=(u\circ\xi)g$ with $\xi\in\cT^1_p(\Omega_1,\Omega_2)$ and non-negative
$g\in\cM(\sW^{1,p}(\Omega_2)\sto\sW^{1,q}(\Omega_2))$ is a vector lattice homomorphism.

The article is organized as follows.
In Section \ref{sec:pre} we fix the setting, give examples and prove preliminary
results. Some of the introduced objects
are well-known,
some are less-known
and some of them are new.
In Section \ref{sec:abstract} we prove the
Abstract Representation Theorem and in Section \ref{sec:sobolev} we apply it
to various Sobolev spaces. In the last and short section 
(Section \ref{sec:examples}) with give some examples.
The sections are split into subsections as follows.
In Subsection \ref{ss:fs} we introduce well-known classes of functions and
in Subsection \ref{ss:cap} we will shortly introduce the classical
$p$-capacity $\Cap_p$ with references to the literature for more informations.
The definitions in Subsection \ref{ss:equi} are new. Here we introduce various equivalence
classes of functions - functions which are not defined everywhere.
It will be important for the Abstract Representation Theorem (Theorem \ref{thm:art}) to
distinguish between pointwise defined functions and equivalence classes of functions.
The relative $p$-capacity is introduced in Subsection \ref{ss:rcap}.
The only use of the relative $p$-capacity is to handle lattice homomorphisms involving
Sobolev spaces with non-vanishing boundary values, such as $\tsW^{1,p}$.
In Subsection \ref{ss:mol} we recall well-known results about the mollification of
$L^p$ and Sobolev functions which we need in Subsection \ref{ss:reg} to deduce that
$L^p$-spaces and Sobolev spaces are regularizable, a notion which is also defined there.
Properties of lattice homomorphisms between Sobolev spaces are given in Subsection \ref{ss:lattice}.
A partition of unity consisting of functions in certain function spaces is introduced
in Subsection \ref{ss:partition}.
In Section \ref{sec:abstract} we prove the Abstract Representation Theorem,
stating that for certain operators $T$ we have a representation of the form $Tu=(u\circ\xi)g$,
which we also call a {\bf CoMu}-representation ({\bf Co}mposition and {\bf Mu}ltiplication) for $T$.
In Section \ref{sec:sobolev} we apply the Abstract Representation Theorem to various 
Sobolev spaces.
Representations of lattice homomorphisms between Sobolev spaces with vanishing boundary
values are considered in Subsection \ref{ss:h-zero} -- Theorem \ref{thm:zero}.
Representations for Sobolev spaces with non-vanishing boundary values are considered
in Subsection \ref{ss:h-nonzero} -- Theorem \ref{thm:nonzero} -- whereas representations
up to the boundary are considered in Subsection \ref{ss:h-nonzero2} -- Theorem \ref{thm:nonzero2}.

%--------------------------------------------------------------------------------------------------
%--------------------------------------------------------------------------------------------------
%--------------------------------------------------------------------------------------------------

\section{Preliminaries and Setting}\label{sec:pre}

In this article 
$\Omega$ always denotes an open and non-empty subset of $\IR^N$,
$(\metric,\dist)$ denotes a metric space,
$\topo$ denotes a topological space and
$\domain$ denotes an arbitrary non-empty set.

%--------------------------------------------------------------------------------------------------
%--------------------------------------------------------------------------------------------------
%--------------------------------------------------------------------------------------------------

\subsection{Function Spaces}\label{ss:fs}

By $C(\topo)$ we denote the space of all real-valued and continuous functions on $\topo$
and by $C_c(\topo)$ the subspace of $C(\topo)$ consisting of those functions having compact
support. By $\cD(\Omega)$ we denote the space of all {\em test functions} on $\Omega$, that is,
\[ \cD(\Omega)  := C^\infty(\Omega)\cap C_c(\Omega) = \klg{u\in C^\infty(\Omega): \supp(u)\subset\Omega \mbox{ is compact}}.
\]
Its topological dual is denoted by $\cD'(\Omega)$ and is called the space of {\em distributions}.
For $p\in [1,\infty)$ the first order Sobolev space $W^{1,p}(\Omega)\subset L^p(\Omega)$
is given by
\begin{eqnarray*}
   W^{1,p}(\Omega) &:=& \klg{\su\in L^p(\Omega): D^\alpha \su\in L^p(\Omega)\mbox{ in }\cD'(\Omega)\;
                             \mbox{ for all }\alpha\in\IN_0^N, \betrag{\alpha}\leq 1},\\
   \norm{\su}^p_{W^{1,p}(\Omega)} &:=& \sum_{\betrag{\alpha}\leq 1} \norm{D^\alpha \su}^p_{L^p(\Omega)}.
\end{eqnarray*}

\subsection{The classical $p$-Capacity}\label{ss:cap}
For $p\in(1,\infty)$ the classical $p$-capacity $\Cap_p$ of a set $A\subset\IR^N$ is given by
\[ \Cap_p(A) := \inf\klg{\norm{\su}^p_{W^{1,p}(\IR^N)}: \su\geq 1\mbox{ a.e. on a neighbourhood of }A}.
\]
A pointwise defined function $u:A\to\IR$ is called
{\em $\Cap_p$-quasi continuous} on $A$ if for each $\varepsilon>0$ there exists an open set
$V\subset\IR^N$ with $\Cap_p(V)\leq\varepsilon$ such that $u$ restricted to $A\setminus V$ is continuous.
A set $P\subset\IR^N$ is called {\em $\Cap_p$-polar} if $\Cap_p(P)=0$ and we say that a property holds
$\Cap_p$-quasi everywhere (briefly $p$-q.e) if it holds except for a $\Cap_p$-polar set.
For more details we refer to
 Adams and Hedberg \cite{adams:96:fsp},
 Biegert \cite{biegert:08:trc}.
 Bouleau and Hirsch \cite{bouleau:91:df},
 Federer and Ziemer \cite{federer:72:tls},
 Fukushima and {\=O}shima and Takeda \cite{fukushima:94:df},
 Mal{\'y} and Ziemer \cite{ziemer:97:fr},
 Maz'ya \cite{mazya:85:ssp},
 Meyers \cite{meyers:75:cpp}
 and the references therein.

\begin{theorem}{Adams and Hedberg \cite[Proposition 6.1.2 and Theorem 6.1.4]{adams:96:fsp} or
                Mal{\'y} and Ziemer \cite[Theorem 2.20 and Corollary 2.23]{ziemer:97:fr}}.
  For every $p\in(1,\infty)$ and $\su\in W^{1,p}(\Omega)$ there exists a $\Cap_p$-quasi continuous
  representative $u$ of $\su$. Such a representative is unique up to a $\Cap_p$-polar set and is denoted by
  $\widetilde\su$.
\end{theorem}

\begin{theorem}{Mal{\'y} and Ziemer \cite[Corollary 2.25]{ziemer:97:fr}}.
  For an arbitrary set $A\subset\IR^N$ and $p\in(1,\infty)$ the $p$-capacity of $A$ is given by
  \[ \Cap_p(A) = \inf\klg{\norm{\su}^p_{W^{1,p}(\IR^N)}:\su\in W^{1,p}(\IR^N),\widetilde\su\geq 1\mbox{ $p$-q.e. on } A}.
  \]
\end{theorem}

\begin{theorem}{Bouleau and Hirsch \cite[Proposition 8.2.5]{bouleau:91:df}}.\label{thm:subseq}
  Let $p\in(1,\infty)$ and $\su_n\in W^{1,p}(\Omega)$ be a sequence which converges in $W^{1,p}(\Omega)$
  to $\su\in W^{1,p}(\Omega)$. Then there exist a $\Cap_p$-polar set $P$ and a subsequence $(\su_{n_k})_k$
  of $(\su_n)_n$ such that $\widetilde \su_{n_k}\to \widetilde \su$ everywhere on $\Omega\setminus P$.
\end{theorem}

%--------------------------------------------------------------------------------------------------
%--------------------------------------------------------------------------------------------------
%--------------------------------------------------------------------------------------------------

\subsection{Equivalence Classes of Functions}\label{ss:equi}

By $\cF(\domain)$ we denote the space of all real-valued functions $f:\domain\to\IR$.
The power set of $\domain$ is denoted by $\cP(\domain)$. We call a subset $\cN\subset\cP(\domain)$
a {\em nullspace} on $\domain$ if it contains the empty set and if it is closed with respect to
countable unions, that is, $\emptyset\in\cN$ and
\[ N_n\in\cN \mbox{ for all } n\in\IN \quad\Longrightarrow\quad \bigcup_{n=1}^\infty N_n\in\cN.
\]
If $\cN$ is a nullspace on $\domain$, then an equivalence relation $\sim_{\cN}$ on $\cF(\domain)$ is given by
\[ f\sim_{\cN} g :\Longleftrightarrow \mbox{ there exists } N\in\cN\mbox{ such that } f=g \mbox{ on } \domain\setminus N.
\]
In the following we consider subspaces $U$ of the quotient space
$\sF(\domain,\cN)$ given by
\[ \sF(\domain,\cN):=\cF(\domain)/\sim_{\cN}.
\]

\begin{definition}
  Let $\cN$ be a nullspace on $\domain$. Then the vector space $\sF(\domain,\cN)$ is equipped with the
  order relation $\leq$ defined by
  \[ \su\leq\sv :\Longleftrightarrow \mbox{ there exist } u\in\su, v\in\sv, N\in\cN
     \mbox{ such that } u\leq v\mbox{ everywhere on }
     \domain\setminus N. 
  \]
\end{definition}

\begin{remark}
  Note that with this ordering, the space $\sF(\domain,\cN)$ is a
  $\sigma$-Dedekind complete vector lattice. For more details we
  we refer to Aliprantis and Burkinshaw \cite[Ch.1,Sect.1]{aliprantis:85:po}.
\end{remark}

\begin{example}
  Let $\cN_0(\Omega)\subset\cP(\Omega)$ denote the set of all nullsets $N\subset\Omega$ with respect
  to the Lebesgue measure. Then $U:=L^p(\Omega)$ is a subspace of $\sF(\Omega,\cN_0(\Omega))$ for every
  $p\in[1,\infty]$.
\end{example}

\begin{example}\label{ex:w1p}
  Let $p\in (1,\infty)$ and let $\cN_p(\Omega)\subset\cP(\Omega)$ consist of all $\Cap_p$-polar sets $N\subset\Omega$.
  Then we define the ({\em refined}) Sobolev spaces $\sW^{1,p}(\Omega),\sW^{1,p}_0(\Omega)\subset \sF(\Omega,\cN_p)$ as follows:
  \begin{eqnarray*} 
      \sW^{1,p}(\Omega) &:=& \klg{[u]_{\cN_p}:u\in\su\in W^{1,p}(\Omega)\mbox{ is $\Cap_p$-quasi continuous}}.\\
      \sW^{1,p}_0(\Omega) &:=& \overline{\cD(\Omega)}^{\sW^{1,p}(\Omega)}.
  \end{eqnarray*}
  Here $[u]_{\cN_p}$ denotes the equivalence class of $u\in\cF(\Omega)$ with respect to $\sim_{\cN_p}$.
\end{example}

\begin{remark}
  Let $\cN$ be a nullspace on $\domain$ and let $U$ be a subspace of $\sF(\domain,\cN)$.
  Saying that a function $u\in\cF(\domain)$ belongs to $U$ means that $[u]_{\cN}$ belongs
  to $U$. For example, by this identification we have $\cD(\Omega)\subset L^p(\Omega)$.
  Moreover, let $\cN_1$ and $\cN_2$ be nullspaces on $\domain$ and let $U_2$ be a subspace of
  $\sF(\domain,\cN_2)$. Saying that $\su\in\sF(\domain,\cN_1)$ belongs to $U_2$ means
  that there exists $u\in\su$ such that $[u]_{\cN_2}\in U_2$. For example, by
  this identification we have $\cD(\Omega)\subset L^p(\Omega)$ and $W^{1,N+\varepsilon}(\Omega)\subset C(\Omega)$
  where $\cD(\Omega)$ and $C(\Omega)$ are identified with a subspace of $\sF(\Omega,\klg{\emptyset})$
  via the previous identification.
\end{remark}

%--------------------------------------------------------------------------------------------------
%--------------------------------------------------------------------------------------------------
%--------------------------------------------------------------------------------------------------

\subsection{The relative $p$-Capacity}\label{ss:rcap}
 In this subsection we introduce the relative $p$-capacity with respect to an open set $\Omega\subset\IR^N$.
 The notion of the relative $2$-capacity was first introduced by Arendt and Warma
 in \cite{arendt:03:lrb} to study the Laplacian with Robin boundary conditions on
 arbitrary domains in $\IR^N$. This notion was extended to $p\in (1,\infty)$ 
 by Biegert in \cite{biegert:08:trc} where also further properties are proved. The importance of the
 relative $p$-capacity is that Sobolev functions in $\tsW^{1,p}(\Omega)$ admit a trace on
 $\partial\Omega$ for {\bf every} open set $\Omega\subset\IR^N$.

\begin{definition}
  For $p\in(1,\infty)$ we let $\tsW^{1,p}(\Omega)$ be the closure of $\sW^{1,p}(\Omega)\cap C_c(\overline\Omega)$
  in $\sW^{1,p}(\Omega)$. Then the relative $p$-capacity $\Cap_{p,\Omega}$ of an arbitrary set $A\subset\overline\Omega$ is given by
  \[ \Cap_{p,\Omega}(A) := \inf\klg{\norm{\su}^p_{\sW^{1,p}(\Omega)}: \su\in\mathcal{Y}(A)}
  \]
  where $\mathcal{Y}(A) := \klg{\su\in\tsW^{1,p}(\Omega):
                            \exists O\subset\IR^N\mbox{ open},\; A\subset O,\; \su\geq 1\mbox{ a.e. on } O\cap\Omega}$.
\end{definition}

\begin{remark}
  Note that in the definition above the intersection is given by
  \[ \sW^{1,p}(\Omega)\cap C_c(\overline\Omega):=
     \klg{u|_\Omega\in\sW^{1,p}(\Omega):u\in C_c(\overline\Omega)}.
  \]
  We should also remark that $\sW^{1,p}_0(\Omega)\subset\tsW^{1,p}(\Omega)$,
  $\tsW^{1,p}(\IR^N)=\sW^{1,p}(\IR^N)$
  and $\Cap_{p,\IR^N}$ is the classical $p$-capacity $\Cap_p$. Moreover,
  if $\Omega$ is a Lipschitz domain or more generally of class $C^0$, then $\tsW^{1,p}(\Omega)$
  and $\sW^{1,p}(\Omega)$ coincide.
\end{remark}

A pointwise defined function $u:\overline\Omega\to\IR$ is called $\Cap_{p,\Omega}$-quasi continuous
if for each $\varepsilon>0$ there exists an open set $V$ in the metric space $\overline\Omega$
with $\Cap_{p,\Omega}(V)\leq\varepsilon$ such that $u$ restricted to $A\setminus V$ is continuous.
A set $P\subset\overline\Omega$ is called $\Cap_{p,\Omega}$-polar if $\Cap_{p,\Omega}(P)=0$ and we
say that a property holds $\Cap_{p,\Omega}$-quasi everywhere (briefly $(p,\Omega)$-q.e.) if it holds
except for a $\Cap_{p,\Omega}$-polar set.

\begin{theorem}Biegert \cite[Theorem 3.22]{biegert:08:trc}.\label{thm:biegert}
  For every $\su\in\tsW^{1,p}(\Omega)$ there exists a $\Cap_{p,\Omega}$-quasi continuous
  function $\tilde u:\overline\Omega\to\IR$ such that $\tilde u=\su$ $\Cap_p$-quasi everywhere on $\Omega$.
  Such a function is unique up to a $\Cap_{p,\Omega}$-polar set.
\end{theorem}

\begin{remark}
  Let $\cN_p^\star(\Omega)$ be the set of all $\Cap_{p,\Omega}$-polar sets in
  $\overline\Omega$. Note that for $A\subset\Omega$
  we have that $\Cap_p(A)=0$ if and only if $\Cap_{p,\Omega}(A)=0$.
  This (together with Theorem \ref{thm:biegert}) shows that we can {\em extend} every
  function $\su\in\tsW^{1,p}(\Omega)$ defined on $\Omega$ in a unique way to a
  $\Cap_{p,\Omega}$-quasi continuous function in $\sF(\overline\Omega,\cN_p^\star(\Omega))$.
  In the following we consider $\tsW^{1,p}(\Omega)$ as a subspace of
  $\sF(\overline\Omega,\cN_p^\star(\Omega))$.
\end{remark}

\begin{theorem}Biegert \cite[Theorem 3.29]{biegert:08:trc}.\label{thm:biegert-cap}
  For an arbitrary set $A\subset\overline\Omega$ and $p\in (1,\infty)$ the relative $p$-capacity of
  $A$ is given by
  \[ \Cap_{p,\Omega}(A) = \inf\klg{\norm{\su}^p_{\tsW^{1,p}(\Omega)}:\su\in\tsW^{1,p}(\Omega),\;\su\geq 1 
      \mbox{ $\Cap_{p,\Omega}$-q.e. on }A}.
  \]
\end{theorem}

\begin{definition}{\bf (Choquet capacity)} Doob\cite[A.II.1]{doob:01:cpt}
  A set function $\sC:\cP(\topo)\to [0,\infty]$ is called a normed {\em Choquet capacity} on $\topo$ if
  it satisfies the following four conditions.
  \begin{itemize}
    \item $\sC(\emptyset)=0$;
    \item $A\subset B\subset\topo$ implies $\sC(A)\leq\sC(B)$;
    \item $(A_n)_n\subset\topo$ increasing implies $\sC(\bigcup_n A_n)=\lim_n \sC(A_n)$;
    \item $(K_n)_n\subset\topo$ decreasing and $K_n$ compact imply $\sC(\bigcap_n K_n)=\lim_n \sC(K_n)$.
  \end{itemize} 
\end{definition}

\begin{theorem}Biegert\cite[Theorem 3.4]{biegert:08:trc}.\label{thm:choquet}
  For an open and non-empty set $\Omega\subset\IR^N$ and $p\in (1,\infty)$ the relative
  $p$-capacity $\Cap_{p,\Omega}$ is a normed Choquet capacity on $\overline\Omega$ and
  \[ \Cap_{p,\Omega}(A) = \inf\klg{\Cap_{p,\Omega}(U):U\mbox{ open in }\overline\Omega\mbox{ and }A\subset U}.
  \]
\end{theorem}

\begin{theorem}Biegert\cite[Proposition 3.5]{biegert:08:trc}.\label{thm:biegert-cpt}
  For a compact set $K\subset\overline\Omega$ and $p\in (1,\infty)$ the relative $p$-capacity
  of $K$ is given by
  \[ \Cap_{p,\Omega}(K) = \inf\klg{\norm{u}^p_{\tsW^{1,p}(\Omega)}:
                      u\in\tsW^{1,p}(\Omega)\cap C_c(\overline\Omega),u\geq 1\mbox{ on }K}.
  \]
\end{theorem}

\begin{theorem}Biegert\cite[Theorem 3.24]{biegert:08:trc}.\label{thm:biegert-new}
  Let $\su_n\in\tsW^{1,p}(\Omega)$ be a sequence which converges in $\tsW^{1,p}(\Omega)$
  to $\su\in\tsW^{1,p}(\Omega)$. Then there exists a subsequence $\su_{n_k}$ which converges
  $\Cap_{p,\Omega}$-quasi everywhere on $\overline\Omega$ to $\su$.
\end{theorem}

%--------------------------------------------------------------------------------------------------
%--------------------------------------------------------------------------------------------------
%--------------------------------------------------------------------------------------------------

\subsection{Mollification}\label{ss:mol}

For $x\in\metric$ and $r>0$ we denote by $B_{\metric}(x,r):=\klg{y\in\metric:\dist(x,y)<r}$ the open
ball in $\metric$ with center $x$ and radius $r$. If no confusion seems likely, we briefly write
$B(x,r)$ instead of $B_\metric(x,r)$. For a set $A\subset\metric$ and $r>0$
we let $B(A,r)$ and $B(A,-r)$ be the open sets given by
\[ B(A,r):=\bigcup_{x\in A} B(x,r),\qquad
   B(A,-r):=\klg{x\in A:\dist(x,A^c)>r}.
\]
We define a sequence of mollifiers as follows: Let $\rho\in\cD(B(0,1))\subset\cD(\IR^N)$
be a non-negative test function such that $\int \rho=1$. Then for $m\in\IN$ and $x\in\IR^N$ we let
\[ \rho_m(x):=m^N\rho(mx).
\]

\begin{theorem}\label{thm:convl2}
  Let $p\in[1,\infty]$, $u\in\su\in L^P(\Omega)$ and let
  $u_n:\IR^N\to\IR$ be given by
  \[ u_n(x):=(\su\star\rho_n)(x) = \int_\Omega \su(y)\rho_n(x-y)\;dy.
  \]
  Then there exists $N\in\cN_0(\Omega)$ such that $u_n(x)\rightarrow u(x)$ for all $x\in\Omega\setminus N$.
\end{theorem}

\begin{proof}
  It is well-known that $u_n(x)\to u(x)$ whenever $x$ is a Lebesgue point for $u$.
  See for instance Ziemer \cite[Theorem 1.6.1(ii)]{ziemer:89:wdf} or 
  Mal{\'y} and Ziemer\cite[Theorem 1.12]{ziemer:97:fr}.
  It is also well-known that almost every $x\in\Omega$ is a Lebesgue point for $u$.
  See Jost \cite[Corollary 19.18]{jost:05:pma} or Mal{\'y} and Ziemer \cite[Theorem 1.24]{ziemer:97:fr}.
\end{proof}

\begin{theorem}\label{thm:convh1}
  Let $p\in(1,\infty)$, $u\in\su\in\sW^{1,p}(\Omega)$ and let
  $u_n:\IR^N\to\IR$ be given by
  \[ u_n(x):=(\su\star\rho_n)(x) = \int_\Omega \su(y)\rho_n(x-y)\;dy.
  \]
  Then there exists $N\in\cN_p(\Omega)$ (i.e. a $\Cap_p$-polar subset of $\Omega$) such that
  $u_n(x)\rightarrow u(x)$ for all $x\in\Omega\setminus N$.
\end{theorem}

\begin{proof}
  This follows as in the proof of Theorem \ref{thm:convl2} with the additional observation that
  $\Cap_p$-quasi every $x\in\Omega$ is a Lebesgue point for $u$ -- see 
  Adams and Hedberg \cite[Theorem 6.2.1]{adams:96:fsp},
  Federer and Ziemer \cite{federer:72:tls} or
  Mal{\'y} and Ziemer \cite[Theorem 2.55]{ziemer:97:fr}.
\end{proof}

\begin{lemma}
  We have that $S_m\in\cL(L^2(\Omega),C(\overline\Omega))$ for all $m\in\IN$ where $S_m$ is given by
  \[ S_m\su:=\su\star\rho_m.
  \]
\end{lemma}

\begin{proof}
  This is Young's inequality stating that for $f\in L^p(\IR^N)$ and
  $g\in L^q(\IR^N)$ with $1\leq p,q\leq\infty$ and $1/p+1/q=1$ we have that
  \[ \norm{f\star g}_\infty \leq \norm{f}_p\norm{g}_q.
  \]
  The continuity (even on $\IR^N$) of $S_m\su$ follows from Mal{\'y} and Ziemer
  \cite[Theorem 1.12(i)]{ziemer:97:fr}.
\end{proof}

%--------------------------------------------------------------------------------------------------
%--------------------------------------------------------------------------------------------------
%--------------------------------------------------------------------------------------------------

\subsection{Regularizable Spaces}\label{ss:reg}
In this subsection we introduce classes of function spaces, called regularizable spaces,
which have the property that functions therein can be in some sense approximated by
smooth functions.

\begin{definition}
  Let $\cN$ be a nullspace on $\topo$ and let $U$ be a subspace of $\sF(\topo,\cN)$.
  We call $U$ {\em regularizable} if there exists a sequence $(S_m)_m$ of linear and positive
  operators $S_m:U\to C(\topo)$ such that the following holds:
  \begin{center}
    For every $u\in\su\in U$ there exists $N\in\cN$ such that $S_m\su(x)\to u(x)$ for all $x\in\topo\setminus N$.
  \end{center}
  In this case we call the sequence $(S_m)_m$ a {\em regularizer sequence} for $U$.
\end{definition}

\begin{example}\label{ex:lpreg}
  The space $U:=L^p(\Omega)\subset\sF(\Omega,\cN_0)$ is regularizable for every $p\in[1,\infty]$.
  A regularizer sequence $(S_m)_m$ for $U$ is given by $S_m\su:=\su\star\rho_m$ --
  Theorem \ref{thm:convl2}.
\end{example}

\begin{example}
  The space $U:=W^{1,p}(\Omega)\subset\sF(\Omega,\cN_0)$ is regularizable for every $p\in[1,\infty)$.
  A regularizer sequence $(S_m)_m$ for $U$ is given by $S_m\su:=\su\star\rho_m$ --
  Theorem \ref{thm:convl2}.
\end{example}

\begin{example}\label{ex:w1preg}
  The space $U:=\sW^{1,p}(\Omega)\subset\sF(\Omega,\cN_p)$ is regularizable for every $p\in(1,\infty)$.
  A regularizer sequence $(S_m)_m$ for $U$ is given by $S_mu:=\su\star\rho_m$ --
  Theorem \ref{thm:convh1}.
\end{example}

\subsection{Lattice Homomorphisms and Local Operators}\label{ss:lattice}

In this subsection we consider in particular lattice homomorphisms between $L^p$-spaces or
Sobolev spaces and we show that they satisfy the conditions in the Abstract
Representation Theorem (Theorem \ref{thm:art}).

\begin{definition}
  Let $\topo$ be a topological space, $\domain_1\subset\topo$ be a dense subset,
  $\cN_j$ be a nullspace on $\domain_j$ and let $U$ be a subspace of $\sF(\domain_1,\cN_1)$.
  Then a linear operator $T:U\to\sF(\domain_2,\cN_2)$ is called {\em $\topo$-local}, if for all
  $u,v\in U\cap C_c(\topo):=\klg{u|_{\domain_1}\in U:u\in C_c(\topo)}$ with disjoint
  support the product $Tu\cdot Tv=0$ in $\sF(\domain_2,\cN_2)$.
\end{definition}

\begin{definition}
  Let $E,F$ be vector lattices. A linear mapping $T:E\to F$ is called a
  {\em vector lattice homomorphism} or briefly {\em lattice homomorphism}
  if $\betrag{Tu}=T\betrag{u}$ for all $u\in E$. If $T$ is in addition
  bijective, then $T$ is called a {\em lattice isomorphism}.
\end{definition}

\begin{lemma}\label{lem:l2oh}
  Let $1\leq p,q\leq\infty$ and let $T$ be a linear mapping from $L^p(\Omega_1)$ into $L^q(\Omega_2)$.
  Then the following are equivalent.
  \begin{enumerate}
    \item $T$ is a lattice homomorphism.
    \item $T$ is a positive and $\Omega_1$-local operator.
    \item $T$ is a continuous, positive and $\Omega_1$-local operator.
  \end{enumerate}
\end{lemma}

\begin{proof}
  {\bf(1)$\Rightarrow$(2)}. The positivity of $T$ is clear. To show that $T$ is local
  let $u,v\in L^p(\Omega_1)\cap C_c(\Omega_1)=C_c(\Omega_1)$. Then
  $0\leq \betrag{Tu}\wedge \betrag{Tv}=T\betrag{u}\wedge T\betrag{v}=T(\betrag{u}\wedge \betrag{v})=0$
  almost everywhere and hence $Tu\cdot Tv=0$ in $\sF(\Omega_2,\cN_0(\Omega_2))$.
  {\bf (2)$\Rightarrow$(3)}. Using that $L^q(\Omega_2)$ is a Banach lattice and that 
  $T$ is positive, we get from Schaefer \cite[Theorem II.5.3]{schaefer:74:blp} that $T$ is continuous.
  {\bf (3)$\Rightarrow$(1)}.
  Let $\su\in L^p(\Omega_1)$. It suffices to show that $T\su^+\wedge T\su^-=0$ (Schaefer \cite[Proposition II.2.5]{schaefer:74:blp}).
  For this, using the continuity and locality of $T$, it is sufficient to show that there exist functions
  $u_n,v_n\in L^p(\Omega_1)\cap C_c(\Omega_1)=C_c(\Omega_1)$ with disjoint support
  (i.e. $\supp(u_n)\cap\supp(v_n)=\emptyset$) such that $u_n\to \su^+$ and $v_n\to\su^-$ in $L^p(\Omega_1)$.
  To find such sequences let $w\in L^p(\Omega_1)\cap C_c(\Omega_1)=C_c(\Omega_1)$ be such
  that $\norm{w-\su}_p\leq 1/n$. Let
  $\delta>0$ be so small that for $u_n:=(w^+-\delta)^+$ and $v_n:=(w^--\delta)^+$
  one has $\norm{w^+-u_n}\leq 1/n$ and $\norm{w^--v_n}\leq 1/n$. It follows that
  $(u_n)_n$ and $(v_n)_n$ are sequences with the desired properties.
\end{proof}

\begin{lemma}\label{lem:h1oh}
  Let $p,q\in(1,\infty)$ and let $T$ be a linear mapping from $\sW^{1,p}_0(\Omega_1)$ into $\sW^{1,q}(\Omega_2)$.
  Then the following are equivalent.
  \begin{enumerate}
    \item $T$ is a lattice homomorphism.
    \item $T$ is a positive and $\Omega_1$-local operator.
    \item $T$ is a continuous, positive and $\Omega_1$-local operator.
  \end{enumerate}
\end{lemma}

\begin{proof}
  {\bf (1)$\Rightarrow$(2)}. Analogous to the proof of Lemma \ref{lem:l2oh}.
  {\bf (2)$\Rightarrow$(3)}. Use Theorem \ref{thm:automatic} instead of Schaefer \cite[Theorem II.5.3]{schaefer:74:blp}
  in the proof of Lemma \ref{lem:l2oh}.
  {\bf (3)$\Rightarrow$(1)}. Analogous to the proof of Lemma \ref{lem:l2oh}.
\end{proof}

\begin{remark}
  Note that Schaefer \cite[Theorem II.5.3]{schaefer:74:blp} cannot be used in the proof above since $\sW^{1,q}(\Omega_2)$
  is not a Banach lattice.
\end{remark}

\begin{theorem}{Arendt \cite[Appendix]{arendt:84:rpo}}\label{thm:automatic}
  Let $E$ and $F$ be Banach spaces and let $E_+$ and $F_+$ be closed cones
  in $E$ and $F$, respectively. Assume that $E_+$ is generating
  (i.e. $E=E_+-E_+$) and $F_+$ is proper (i.e. $F_+\cap (-F_+)=\klg{0}$).
  If $T:E\to F$ is a linear and positive operator, then $T$ is continuous.
\end{theorem}

\begin{lemma}\label{lem:h1ohp}
  Let $p,q\in (1,\infty)$ and let $T:\tsW^{1,p}(\Omega_1)\to\sW^{1,q}(\Omega_2)$ be linear.
  Then the following are equivalent.
  \begin{enumerate}
    \item $T$ is a lattice homomorphism.
    \item $T$ is a positive and $\overline\Omega_1$-local operator.
    \item $T$ is a continuous, positive and $\overline\Omega_1$-local operator.
  \end{enumerate}
\end{lemma}

\begin{proof}
  {\bf (1)$\Rightarrow$(2)}. Analogous to the proof of Lemma \ref{lem:l2oh}.
  {\bf (2)$\Rightarrow$(3)}. Use Theorem \ref{thm:automatic} instead of Schaefer \cite[Theorem II.5.3]{schaefer:74:blp}
  in the proof of Lemma \ref{lem:l2oh}.
  {\bf (3)$\Rightarrow$(1)}. Analogous to the proof of Lemma \ref{lem:l2oh}.
\end{proof}

%----------------------------------------------------------------------------------------------------------
%----------------------------------------------------------------------------------------------------------
%----------------------------------------------------------------------------------------------------------

\subsection{Partition of Unity}\label{ss:partition}

In this subsection we consider topological spaces which admit a ``Partition of Unity''
of certain function classes. See also Albeverio and Ma and R\"ockner \cite[Definition 1.2]{roeckner:97:pof}.

\begin{definition}
  Let $\topo$ be a 
  %locally compact and $\sigma$-compact 
  topological space
  and let $U$ be a subspace of $C_c(\topo)$.
  Then we say that {\em $\topo$ admits a partition of unity of class $U$}
  or that {\em $U$ is rich in $\topo$} if for any countable open covering $\bigcup_{n\in\IN} O_n=\topo$
  with relatively compact sets $O_n$ there exist functions
  $\varphi_n\in U$ such that
  $0\leq\varphi_n\leq 1$, $\supp(\varphi_n)\subset O_n$, $(\supp(\varphi_n))_n$ is locally finite and 
  $\sum_n \varphi_n(x)=1$ for all $x\in\topo$.
  Note that this definition has only a real meaning if $\topo$ is a locally compact
  and $\sigma$-compact topological space. In fact, otherwise there exist no such open
  covering as considered above and hence every subspace $U$ of $C_c(\topo)$ is rich in $\topo$.
\end{definition}

\begin{remark}
  A collection of sets in a topological space $\topo$ is called {\em locally finite}
  if every $x\in\topo$ has a neighbourhood which intersects only finitely many elements of this collection. 
  See Munkres \cite[Definition 39]{munkres:00:top}.
\end{remark}

\begin{example}
  Let $\metric$ be a locally compact separable metric space
  (hence $\metric$ is a paracompact $\sigma$-compact Hausdorff space).
  Then $\metric$ admits a partition of unity of class $C_c(\metric)$.
  See Munkres \cite[Theorem 41.7]{munkres:00:top}.
\end{example}

\begin{lemma}\label{lem:testricha}
  Let $\Omega\subset\IR^N$ be a non-empty open set and let $\metric\subset\Omega$ be a metric
  space equipped with the metric from $\IR^N$. Then $R:=\klg{u|_\metric:u\in\cD(\Omega)}\cap C_c(\metric)$
  is rich in $\metric$.
\end{lemma}

\begin{proof}
  Let $O_n$ be open and relatively compact sets in $\metric$ such that $\bigcup_n O_n=\metric$.
  Then there exist open sets $U_n$ in $\Omega$ such that $O_n=U_n\cap\metric$.
  For $U:=\bigcup_n U_n$ there exist open sets $W_n$
  and $V_n$ in $U$ such that $\overline{W_n}^U\subset V_n$, $\overline{V_n}^U\subset U_n$,
  $\bigcup_n W_n=U$ and $(V_n)_n$ is locally finite in $U$. See Munkres \cite[Lemma 41.6]{munkres:00:top}.
  Since $K_n:=\overline{W_n\cap\metric}^\metric\subset \overline{W_n}^U \subset V_n$ is compact there exist
  $\psi_n\in\cD(V_n)\subset\cD(U)$ such that $0\leq\psi_n\leq 1$ in $U$ and
  $\psi_n\equiv 1$ on $K_n$. Hence $\supp_\metric(\psi_n|_\metric)\subset O_n$.
  Define $\Psi(x):=\sum_n \psi_n(x)$. Then $\Psi\in C^\infty(U)$ and
  $\Psi\geq 1$ on $\metric$. For $\varphi_n:=(\psi_n/\Psi)|_\metric$ we have the desired properties.
\end{proof}

\begin{corollary}\label{cor:testrich0}
  Let $\Omega\subset\IR^N$ be a non-empty open set. Then $\cD(\Omega)$ is rich in $\Omega$.
\end{corollary}

\begin{proof}
  Apply Lemma \ref{lem:testricha} with $\metric=\Omega$.
\end{proof}

\begin{corollary}\label{cor:testrich1}
  Let $V\subset\IR^N$ be a non-empty open set.
  Then $R:=\klg{u|_{\overline V}:u\in\cD(\IR^N)}\cap C_c(\overline V)$
  is rich in $\overline V$. In particular, $\tsW^{1,p}(V)\cap C_c(\overline V)$
  is rich in $\overline V$.
\end{corollary}

\begin{proof}
  Apply Lemma \ref{lem:testricha} with $\metric:=\overline V\subset\Omega:=\IR^N$.
\end{proof}

For the proof of the Abstract Representation Theorem we need the following.

\begin{lemma}{\bf (Partition of Unity)}\label{lem:partunity}
  Let $\metric$ be a locally compact separable metric space
  and let $U\subset C_c(\metric)$ be a rich subspace in $\metric$.
  Then for $k,n\in\IN$ there exist functions $\varphi_{k,n}\in U$
  satisfying the following properties:
  \begin{itemize}
    \item For $k,n\in\IN$ there exist $z_{k,n}\in\metric$, $r_{k,n}\in(0,1/n]$ such that
          $\varphi_{k,n}\in C_c(B(z_{k,n},r_{k,n}))$.
    \item For $n\in\IN$ and $K\subset\metric$ compact the set
          $\klg{k\in\IN:\supp(\varphi_{k,n})\cap K\not=\emptyset}$ is finite.
    \item For $n\in\IN$ fixed we have (pointwise) that $\sum_{k=1}^\infty \varphi_{k,n}\equiv 1$ on $\metric$.
  \end{itemize}
\end{lemma}

\begin{proof}
  Let $n\in\IN$ be fixed and let $K_j\subset\metric$ be a sequence of compact sets
  whose union is $\metric$. Since $\metric$ is locally compact and $K_j$ is compact,
  we know that $K_j$ is in a finite union of open and relatively compact balls with center
  in $K_j$ and radius in $(0,1/n]$. Since the union of all $K_j$ is $\metric$, we get that
  the countable union of all such open and relatively compact balls, denoted by $B(z_{k,n},r_{k,n})$ for
  $k\in\IN$, equals $\metric$.
  Since $\metric$ admits a partition of unity of class $U$,
  we get a family of functions $(\varphi_{k,n})_k$ in $U$
  satisfying the desired properties.
\end{proof}

%--------------------------------------------------------------------------------------------------
%--------------------------------------------------------------------------------------------------
%--------------------------------------------------------------------------------------------------

\section{Abstract Representation Theorem}\label{sec:abstract}

The aim of this section is to proof a representation for linear, positive and
local operators defined on a rich subspace similar to the following well-known
Theorem --
see Aliprantis and Burkinshaw \cite[Theorem 7.22]{aliprantis:85:po}).

\begin{theorem}
  Let $X$ and $Y$ be two compact Hausdorff spaces.
  Then for a nonzero positive operator $T:C(X)\to C(Y)$ the following
  statements are equivalent:
  \begin{itemize}
    \item $T$ is a lattice homomorphism.
    \item There exist a unique non-negative function $g\in C(Y)$ and a function
          $\xi:Y\to X$ which is continuous on $\klg{y\in Y:g(y)>0}$, such that
          for all $y\in Y$ and all $u\in C(X)$
          \[ Tu(y)=u(\xi(y))g(y).
          \]
  \end{itemize}
\end{theorem}

\begin{definition}
  Let $\domain_j$ be a non-empty set and $\cN_j$ be a nullspace on $\domain_j$ for $j=1,2$ and
  let $T$ be a linear mapping from a subspace $U\subset\sF(\domain_1,\cN_1)$ into
  $\sF(\domain_2,\cN_2)$. Then a pair of functions $(\xi,g)$ where $\xi:\domain_2\to\domain_1$ and
  $g:\domain_2\to\IR$ is called a {\em CoMu-Representation} of $T$ if
  for every $u\in\su\in U$ and every $f\in T\su$ there exists $N\in\cN_2$ such that
  $f(y) = u(\xi(y))g(y)$ for all $y\in\domain_2\setminus N$.
  In this case we briefly write $T\su=(\su\circ \xi)g$.
\end{definition}

\begin{remark}
  Let $u\in\su\in U$ and $N\in\cN_1$ be given. Define $v\in\su$ by $v(x):=u(x)$ if
  $x\in\domain_1\setminus N$ and $v(x):=u(x)+1$ if $x\in N$. Then for a CoMu-Representation
  $(\xi,g)$ of $T$ we have that there exists $Q\in\cN_2$ such that
  $u(\xi(y))g(y) = v(\xi(y))g(y)$ for all $y\in\domain_2\setminus Q$.
  Therefore, on $P:=\klg{y\in\domain_2:g(y)\not=0}\setminus Q$,
  we get that $u(\xi(y))=v(\xi(y))$. This implies that $\xi(y)\not\in N$ for all $y\in P$.
  Hence $\xi^{-1}(N)\subset Q\cup\klg{y\in\domain_2:g(y)=0}$. This shows that
  for all nullsets $N\in\cN_1$ there is a nullset $Q\in\cN_2$ such that
  $\xi^{-1}(N)\cap\klg{y\in\domain_2:g(y)\not=0}\subset Q$.
\end{remark}

\begin{lemma}{\bf (Representation for linear and positive functionals)}\label{lem:riesz}
  Let $\metric$ be a locally compact separable metric space
  and let $U\subset C_c(\metric)$ be a rich subspace in $\metric$.
  If $T:U\to\IR$ is linear and positive, then there exists a Radon measure $\mu$ on
  $\metric$ such that
  \[ Tu = \int_{\metric} u\;d\mu\quad\mbox{ for all } u\in U.
  \]
\end{lemma}

\begin{proof}
  For $v\in C_c(\metric)$ we let $R(v):=\klg{w\in U:v\leq w}$ and
  $p(v):=\inf_{w\in R(v)} Tw$. Then $p(v)=Tv$ for all $v\in U$. We show that $p:C_c(\metric)\to\IR$ is
  a sublinear functional, that is, $p(u+v)\leq p(u)+p(v)$ and $p(\lambda u)=\lambda p(u)$ for all $u,v\in C_c(\metric)$ 
  and all $\lambda\geq 0$. Let $v\in C_c(\metric)$. Since $U$ is rich there exists $\varphi\in U\cap C_c(\metric)$,
  $\varphi\geq 0$ such that $\varphi\equiv 1$ on $\supp(v)$. Hence $-\varphi\norm{v}_\infty \leq v\leq \varphi\norm{v}_\infty$.
  This shows that $p(v)\in\IR$ for all $v\in C_c(\metric)$. For $u,v\in C_c(\metric)$ there exists
  $u_\varepsilon\in R(u)$ and $v_\varepsilon\in R(v)$ such that $Tu_\varepsilon\leq p(u)+\varepsilon$ 
  and $Tv_\varepsilon\leq p(v)+\varepsilon$. Hence $p(u+v)\leq Tu_\varepsilon+Tv_\varepsilon\leq p(u)+p(v)+2\varepsilon$.
  For $\varepsilon\to 0+$ we get $p(u+v)\leq p(u)+p(v)$. Now let
  $\lambda>0$. Then $p(\lambda u)=\inf_{w\in R(\lambda u)} Tw = \inf_{w\in R(u)} \lambda Tw = \lambda p(u)$.
  It follows from the Hahn-Banach Theorem (see Conway \cite[Theorem 6.2]{conway:90:cif}) that there exists
  $\widetilde T:C_c(\metric)\to\IR$ such that $\widetilde Tu=Tu$ for all $u\in U$ and
  $\widetilde Tu\leq p(u)$ for all $u\in C_c(\metric)$. Since $p(u)\leq 0$ for all $u\in C_c(\metric)$, $u\leq 0$
  we get that $\widetilde T$ is positive. By the Riesz-Markov Representation Theorem
  (see Royden \cite[Theorem 13.4.23]{royden:88:ran}) we get that there exists a Radon measure $\mu$ on $\metric$ such that
  $\widetilde Tu=\int_{\metric} u\;d\mu$ for all $u\in C_c(\metric)$. In particular, we get that
  $Tu=\widetilde Tu=\int_{\metric} u\;d\mu$ for all $u\in U$.
\end{proof}

\begin{theorem}{\bf (Abstract Representation Theorem)}\label{thm:art}
  Assume the following.
  \begin{enumerate}
    \item $\metric$ a locally compact separable metric space and
          $U_1\subset C_c(\metric)$ rich in $\metric$.
    \item $\topo$ a topological space, $\cN_\topo$ a nullspace on $\topo$ and
          $U_2\subset\sF(\topo,\cN_\topo)$ regularizable. 
    \item $T:U_1\to U_2$ a linear, positive and $\metric$-local operator.
  \end{enumerate}
  Then $T$ has a CoMu-Representation $(\xi,g)$ with
  $\xi:\topo\to\metric$ and $g:\topo\to [0,\infty)$, that is,
  \[ Tu = (u\circ\xi)g\;\mbox{ in }\sF(\topo,\cN_\topo)
  \]
  for all $u\in U_1$.
\end{theorem}

\begin{proof}
  Let $S_m:U_2\to C(\topo)$ be a regularizer sequence for $U_2$,
  $\varphi_{k,n}\in U_1$ be given from Lemma \ref{lem:partunity} and let
  $\psi_{k,n}\in T\varphi_{k,n}$ be fixed.
  For $T_m:=S_m\circ T$ we get by our assumptions that there exist $N_{k,n}\in\cN_2$ such that
  $T_m\varphi_{n,k}\to \psi_{k,n}$ everywhere on $\topo\setminus N_{k,n}$.
  If $\supp(\varphi_{k,n})\cap\supp(\varphi_{j,m})\not=\emptyset$ we let
  $N_{k,n,j,m}:=\emptyset$, otherwise ($T$ is local) we let $N_{k,n,j,m}\in\cN_2$ be such that
  $\psi_{k,n}(y)\cdot\psi_{j,m}(y)=0$ for all $y\in\topo\setminus N_{k,n,j,m}$.
  Now let $N\in\cN_2$ be the union of all $N_{k,n}$ and $N_{k,n,j,m}$ and let
  \[ \topo':=\klg{y\in\topo\setminus N:\mbox{ there exist } k,n\in\IN\mbox{ such that } \psi_{k,n}(y)>0}.
  \]
  {\flushleft{\bf Step 1:} Radon measures $\mu_{y,m}$.}
  Let $y\in\topo$ and $m\in\IN$ be fixed. By the properties of $T_m$ it follows that
  $\delta_y\circ T_m:U_1\to\IR$ is linear and positive.
  It follows from Lemma \ref{lem:riesz} that there exists a Radon measure $\mu_{y,m}$ on $\metric$ such that
  \[ \klr{T_m u}(y) = \int_\metric u\;d\mu_{y,m}\qquad\mbox{ for all } u\in U_1.
  \]
  {\flushleft\bf Step 2:} We show that for every $y\in\topo'$ there exists
  $\xi(y)\in\metric$ such that for all compact sets $K\subset\metric\setminus\klg{\xi(y)}$ we
  have that
  \[ \mu_{y,m}(K)\to 0\quad\mbox{ as } m\to\infty.
  \]
  For the proof let $y\in\topo'$ be fixed. Then there exist $k_0,n_0\in\IN$ such that
  $\psi_{k_0,n_0}(y)>0$.
  {\flushleft\bf Step 2a}: We show that for every $n\geq n_0$ there exists $k_0(n)\in\IN$ such that $\psi_{k_0(n),n}(y)>0$.
  Let $n\geq n_0$ be fixed. Since $\supp(\varphi_{k_0,n_0})\subset\metric$ is compact, there
  exists $j\in\IN$ such that
  \[ \Phi:=\sum_{k=1}^j \varphi_{k,n}\equiv 1\mbox{ on }\supp(\varphi_{k_0,n_0}).
  \]
  Using that $T_m:U_1\to C(T)$ is positive, we get that
  \[ \psi_{k_0,n_0}(y) \stackrel{m}{\leftarrow} (T_m \varphi_{k_0,n_0})(y) \leq (T_m\Phi)(y) =
     \sum_{k=1}^j (T_m\varphi_{k,n})(y) \stackrel{m}{\rightarrow} \sum_{k=1}^j \psi_{k,n}(y).
  \]
  Therefore there exists $k_0(n)\in\klg{1,\dots,j}$ such that $\psi_{k_0(n),n}(y)>0$.
  Let $z_n:=z_{k_0(n),n}$ denote the center and $r_n:=r_{k_0(n),n}\leq 1/n$ denote the radius
  of the ball $B$ containing the support of $\varphi_{k_0(n),n}$ (see Lemma \ref{lem:partunity}).
  {\flushleft\bf Step 2b}: Let $K\subset\metric$ be a compact set and let
  $K_n:=K\setminus B(z_n,3/n)$ for $n\geq n_0$. We show that $\mu_{y,m}(K_n)\to 0$ as
  $m\to\infty$ for all $n\geq n_0$.
  Let $n\geq n_0$ be fixed. Then there exists $j\in\IN$ such that
  \[ \Phi:=\sum_{k=1}^j \varphi_{k,n} \equiv 1\mbox{ on } K.
  \]
  For $M_n:=\klg{k\in\IN:1\leq k\leq j, \supp(\varphi_{k,n})\cap\supp(\varphi_{k_0(n),n})=\emptyset}$ we get that
  \[ \eta_n:=\sum_{k\in M_n} \varphi_{k,n} \equiv 1\mbox{ on } K_n
     \quad\mbox{and}\quad
     (\delta_y\circ T_m)\eta_n = \sum_{k\in M_n} (\delta_y\circ T_m)\varphi_{k,n} \stackrel{m}{\to} \sum_{k\in M_n} \psi_{k,n}(y).
  \]
  Since $\supp(\varphi_{k,n})\cap\supp(\varphi_{k_0(n),n})=\emptyset$ for all $k\in M_n$, we get that
  $\psi_{k,n}(y)\psi_{k_0(n),n}(y)=0$ and since $\psi_{k_0(n),n}(y)>0$ it follows that
  $\psi_{k,n}(y)=0$ for all $k\in M_n$. Therefore
  \[ \mu_{y,m}(K_n) \leq \int \eta_n\;d\mu_{y,m} = (\delta_y\circ T_m)\eta_n \rightarrow
     \sum_{k\in M_n} \psi_{k,n}(y) = 0.
  \]

  {\flushleft\bf Step 2c}: We show that $(z_n)_n$ is a Cauchy sequence in $\metric$.
  Let $\varepsilon>0$ be fixed and let $m_0\geq n_0$ be such that $m_0\geq 6/\varepsilon$.
  Assume that there exist $n,l\geq m_0$ such that $\dist\klrk{z_l,z_n}\geq \varepsilon\geq 6/m_0$,
  then
  \[ B(z_l,3/l)\cap B(z_n,3/n)\subset
     B(z_l,3/m_0)\cap B(z_n,3/m_0)=\emptyset.
  \]
  Let $K:=\supp(\varphi_{k_0,n_0})$. Then
  \[ 0<\psi_{k_0,n_0}(y) = \lim_m \int \varphi_{k_0,n_0}\;d\mu_{y,m}
     \leq \limsup_m \mu_{y,m}(K) \leq \limsup_m \mu_{y,m}(K_l)+\mu_{y,m}(K_n)=0,
  \]
  a contradiction. Let $\xi(y):=\lim_n z_n\in\overline\metric$ where $(\overline\metric,\overline\dist)$
  denotes a completion of $(\metric,\dist)$.
 
  {\flushleft\bf Step 2d}: We show the assertion of Step 2.
  Let $K\subset\metric\setminus\klg{\xi(y)}$ be a compact set and $\delta:=\overline\dist(K,\xi(y))>0$.
  Let $n\geq \max(n_0,6/\delta)$ be such that $z_n\in B(\xi(y),\delta/2)$. Then
  \[ B(z_n,3/n)\cap K=\emptyset \quad\text{ whence }\quad K_n:=K\setminus B(z_n,3/n)=K.
  \]
  Hence by Step 2b we get that $\mu_{y,m}(K)\to 0$ as $m\to\infty$.
  Assume that $\xi(y)\in \overline\metric\setminus\metric$. Then
  \[ 0<\psi_{k_0,n_0}(y) = \lim_m \int \varphi_{k_0,n_0}\;d\mu_{y,m} \leq
    \lim_m \mu_{y,m}(\supp(\varphi_{k_0,n_0}))=0,
  \]
  a contradiction and hence $\xi(y)\in\metric$.

  {\flushleft\bf Step 3}: We show that for $y\in\topo'$ there exists $g(y)\in (0,\infty)$
  such that $\int_\metric u d\mu_{y,m}\to g(y)u(\xi(y))$ for all $u\in U_1$. For this
  let $\omega\subset\metric$ be an open and relatively compact set containing $\xi(y)$.
  Then there exist $j\in\IN$ such that
  \[ \Phi:=\sum_{k=1}^j \varphi_{k,n_0} \equiv 1\mbox{ on } \overline\omega.
  \]
  It follows that 
  \begin{align*}
     g(y):=\lim_m \mu_{y,m}(\omega) &= \lim_m \int_{\omega} \Phi\;d\mu_{y,m} 
       = \lim_m \int_{\supp(\Phi)} \Phi\;d\mu_{y,m} - \int_{\supp(\Phi)\setminus\omega} \Phi\;d\mu_{y,m} \\
       &= \lim_m \int_{\metric} \Phi\;d\mu_{y,m} = \sum_{k=1}^j \psi_{k,n_0}(y)\in [0,\infty).
  \end{align*}
  Note that $g(y)$ does not depend on the particular choice of $\omega$.
  Assume that $g(y)=0$, then $0<\psi_{k_0,n_0}(y)=\lim_m \int_\omega \varphi_{k_0,n_0}\;d\mu_{y,m} =0$,
  a contradiction and hence $g(y)>0$.
  Now let $u\in U_1$ and $\varepsilon>0$. By the continuity of $u$ there exists
  $\delta>0$ such that $\betrag{u(\xi(y))-u(x)}\leq\varepsilon$ for all
  $x\in \omega:=B(\xi(y),\delta)$. Without loss of generality we may assume that
  $\omega$ is relatively compact. Hence
  \begin{eqnarray*}
     \limsup_m \int u\;d\mu_{y,m} &=& 
        \limsup_m \int_\omega u\;d\mu_{y,m}
        \leq \limsup_m \mu_{y,m}(\omega)[u(\xi(y))+\varepsilon] \\
     &=& g(y)[u(\xi(y))+\varepsilon] = \liminf_m \mu_{y,m}(\omega)[u(\xi(y))+\varepsilon] \\
     &\leq& \liminf_m \int_\omega u + 2\varepsilon\;d\mu_{y,m} =
             2\varepsilon g(y) + \liminf_m \int u\;d\mu_{y,m}.
  \end{eqnarray*}
  Since $\varepsilon>0$ was arbitrary, the claim follows.

  {\flushleft\bf Step 4}: We finish the proof of the theorem.
  Let $u\in U_1$ and $f\in Tu$ be fixed. Then there exists $\hat N\in\cN_2$
  such that $T_m u\to f$ everywhere on $\topo\setminus\hat N$. Let $M:=N\cup\hat N$
  and $x_0\in\metric$ be fixed. For $y\in\topo\setminus\topo'$ we let $\xi(y):=x_0$
  and $g(y):=0$. We show that for all $y\in\topo\setminus M$ we have that
  \[ f(y) = u(\xi(y))g(y).
  \]
  Let $y\in\topo\setminus M$ be fixed. Then there are two possibilities,
  $y\in\topo'$ or $y\not\in\topo'$. If $y\not\in\topo'$ then for all
  $k,n\in\IN$ we have that $\psi_{k,n}(y)=0$. Let $j\in\IN$ be such that
  \[ \Phi:=\sum_{k=1}^j \varphi_{k,1} \equiv 1\mbox{ on } \supp(u).
  \]
  Then
  \[ \betrag{f(y)} = \lim_m \betrag{T_m u(y)} \leq
     \lim_m \int \betrag{u}\;d\mu_{y,m} \leq
     \norm{u}_\infty \lim_m \int \Phi\;d\mu_{y,m} =\norm{u}_\infty \sum_{k=1}^j \psi_{k,1}(y) = 0.
  \]
  Therefore (since $g(y)=0$) we get that $f(y)=0=u(\xi(y))g(y)$. If $y\in\topo'$ then
  \[ f(y)=\lim_m (T_mu)(y) = \lim_m \int u\;d\mu_{y,m} = u(\xi(y))g(y).
  \]
\end{proof}

%--------------------------------------------------------------------------------------------------
%--------------------------------------------------------------------------------------------------
%--------------------------------------------------------------------------------------------------
%--------------------------------------------------------------------------------------------------
%--------------------------------------------------------------------------------------------------

\section{Representation of Lattice Homomorphisms}\label{sec:sobolev}

In this section we apply the Abstract Representation Theorem (Theorem \ref{thm:art}) to lattice
homomorphisms between $L^p$ and Sobolev spaces defined on open and non-empty sets $\Omega$ in $\IR^N$.
This was the main motivation for the work we did in the previous section.

%--------------------------------------------------------------------------------------------------
%--------------------------------------------------------------------------------------------------
%--------------------------------------------------------------------------------------------------

\subsection{Sobolev Spaces with Vanishing Boundary Values}\label{ss:h-zero}

Let $\Omega_1,\Omega_2\subset\IR^N$ be non-empty open sets and let $p,q\in(1,\infty)$.
In this subsection we assume that $T:\sW^{1,p}_0(\Omega_1)\to\sW^{1,q}(\Omega_2)$
is a lattice homomorphism. It follows from
\begin{itemize}
  \item Example \ref{ex:w1preg} that $\sW^{1,q}(\Omega_2)$ is regularizable,
  \item Lemma \ref{lem:h1oh} that $T$ is $\Omega_1$-local and positive (and continuous),
  \item Corollary \ref{cor:testrich0} that $\sW^{1,p}_0(\Omega_1)\cap C_c(\Omega_1)$ is rich in $\Omega_1$ and
  \item Theorem \ref{thm:art} that $T|_{\sW^{1,p}_0(\Omega_1)\cap C_c(\Omega_1)}$ has a CoMu-Representation $(\xi,g)$,
\end{itemize}
that is, for all $u\in\sW^{1,p}_0(\Omega_1)\cap C_c(\Omega_1)\supset\cD(\Omega_1)$ and $f\in Tu$
there exists a $\Cap_q$-polar set $N$ such that
\begin{equation}\label{eq:ident}
  f(y) = u(\xi(y))g(y)\quad\mbox{ for all } y\in\Omega_2\setminus N.
\end{equation}
In order to extend Equation \eqref{eq:ident} to $u\in\su\in\sW^{1,p}_0(\Omega_1)$ we need the
following lemmata.

\begin{lemma}\label{lem:cap}
  Let $K_j\subset\Omega_j$ be compact sets and let $G_m:=\klg{y\in\Omega_2:g(y)>1/m}$ for $m\in\IN$.
  Then there exists a constant $C=C(K_1,K_2)$ such that for every compact set $K\subset K_1$ the
  following estimate holds:
  \[ \Cap_q(\xi^{-1}(K)\cap G_m\cap K_2)\leq C^q m^q\norm{T}^q\Cap_p(K)^{q/p}.
  \]
\end{lemma}

\begin{proof}
  Let $\psi_j\in\cD(\Omega_j)$ be such that $\psi_j\geq 1$ on $K_j$ and let
  $\varphi_n\in\cD(\IR^N)$ be such that $\varphi_n\geq 1$ on $K$ and
  $\norm{\varphi_n}^p_{W^{1,p}(\IR^N)}\leq \Cap_p(K)+1/n$.
  Let $f_n\in T(\varphi_n\psi_1)$ be fixed and let $N$ be a $\Cap_q$-polar set such that
  \[ f_n(y) = (\varphi_n\psi_1)(\xi(y))g(y)\qquad\mbox{ for all } y\in\Omega_2\setminus N,\;n\in\IN.
  \]
  Then for $y\in\xi^{-1}(K)\cap G_m\cap K_2\cap N^c$ we get that
  \[ m\psi_2(y)f_n(y) \geq mg(y) \geq 1.
  \]
  Hence (using that $m\psi_2f_n$ is $\Cap_q$-quasi continuous) we get that
  \begin{eqnarray*}
     \Cap_q(\xi^{-1}(K)\cap G_m\cap K_2) 
        &\leq& \norm{mf_n\psi_2}^q_{W^{1,q}(\Omega_2)} \leq
     m^q C_2^q \norm{T}^q\norm{\varphi_n\psi_1}^q_{W^{1,p}(\Omega_1)} \\
        &\leq& m^qC_2^qC_1^q\norm{T}^q\norm{\varphi_n}_{W^{1,p}(\IR^N)}^q.
  \end{eqnarray*}
  For $n\to\infty$ the claim follows.
\end{proof}

\begin{lemma}
  Let $K_j\subset\Omega_j$ be compact sets and let $G_m:=\klg{y\in\Omega:g(y)>1/m}$
  for $m\in\IN$. Then there exists a constant $C=C(K_1,K_2)$ such that for all sets
  $M$ in the interior $K_1^\circ$ of $K_1$ the following estimate holds:
  \[ \Cap_q(\xi^{-1}(M)\cap G_m\cap K_2)\leq C^q m^q\norm{T}^q\Cap_p(M)^{q/p}.
  \]
\end{lemma}

\begin{proof}
  Let $\varepsilon>0$. Then there exists an open set $O\subset K_1^\circ$ containing
  $M$ such that $\Cap_p(O)\leq\Cap_p(M)+\varepsilon$. Let $C_n\subset O$ be an increasing sequence
  of compact sets such that $\bigcup_n C_n=O$. By Lemma \ref{lem:cap} we get that
  \begin{eqnarray*}
    \Cap_q(\xi^{-1}(M)\cap G_m\cap K_2) &\leq&
    \Cap_q(\xi^{-1}(O)\cap G_m\cap K_2) \\
      &=& \lim_n \Cap_q(\xi^{-1}(C_n)\cap G_m\cap K_2) \\
      &\leq& \lim_n C^qm^q\norm{T}^q\Cap_p(C_n)^{q/p} \\
      &=& C^qm^q\norm{T}^q\Cap_p(O)^{q/p} \\
      &\leq& C^qm^q\norm{T}^q[\Cap_p(M)+\varepsilon]^{q/p}.
  \end{eqnarray*}
  For $\varepsilon\to 0+$ the claim follows.
\end{proof}

\begin{lemma}\label{lem:polar}
  The set $\xi^{-1}(P)\cap\klg{y\in\Omega_2:g(y)>0}$ is $\Cap_q$-polar for every $\Cap_p$-polar
  set $P\subset\Omega_1$. 
\end{lemma}

\begin{proof}
  Let $\omega^{j}_n\subset\Omega_j$ be increasing sequences of bounded open sets such that
  $\overline\omega_n^{j}\subset\omega_{n+1}^{j}$ and
  $\bigcup_n \omega^{j}_n = \Omega_j$. For all $n,k,m\in\IN$ we get by the previous lemma that
  \[ \Cap_q(\xi^{-1}(P\cap \omega^{1}_n)\cap G_m\cap \omega_k^{2}) \leq C^q_{n,k}m^q\norm{T}^q\Cap_p(P\cap \omega^{1}_n)=0.
  \]
  Now the claim follows from the identity
  \[ \xi^{-1}(P)\cap\klg{y\in\Omega_2:g(y)>0}=\bigcup_{n,m,k} \xi^{-1}(P\cap\omega^{1}_n)\cap G_m\cap \omega_k^{2}.
  \]
\end{proof}

The following theorem is one of the main theorems in this article.
It says that every lattice homomorphism between Sobolev spaces admits a CoMu-Representation.

\begin{theorem}\label{thm:zero}
  Let $\Omega_1,\Omega_2\subset\IR^N$ be non-empty open sets and let $p,q\in(1,\infty)$.
  Assume that $T:\sW^{1,p}_0(\Omega_1)\to\sW^{1,q}(\Omega_2)$ is a lattice homomorphism.
  Then there exists a CoMu-Representation $(\xi,g)$ of $T$ with $\xi:\Omega_1\to\Omega_2$ and $g:\Omega_2\to [0,\infty)$,
  that is,
  \[ T\su = (\su\circ\xi)g\qquad\mbox{ for all } \su\in\sW^{1,p}_0(\Omega_1).
  \]
  More precisely, this means that for every $u\in\su\in\sW^{1,p}_0(\Omega_1)$ and every $f\in T\su$ there exists
  a $\Cap_q$-polar set $N\subset\Omega_2$ such that
  \[ f(y) = u(\xi(y))g(y)\quad\mbox{ for all } y\in\Omega_2\setminus N.
  \]
\end{theorem}

\begin{proof}
  Let $u_n\in\cD(\Omega_1)$ be a sequence of test functions converging in
  $\sW^{1,p}(\Omega_1)$ to $\su\in\sW^{1,p}_0(\Omega_1)$. For $u\in\su$
  (after passing to a subsequence, Theorem \ref{thm:subseq})
  there exists a $\Cap_p$-polar set $P$ such that 
  $u_n\to u$ everywhere on $\Omega_1\setminus P$. Now let $f_n\in Tu_n$ and $f\in T\su$ be fixed.
  Then (after passing to a subsequence, Theorem \ref{thm:subseq})
  there exists a $\Cap_q$-polar set $N_1$ such that $f_n\to f$ everywhere on
  $\Omega_2\setminus N_1$. Let $N_2$ be a $\Cap_q$-polar set such that the
  following holds for all $y\in\Omega_2\setminus N_2$ and all $n\in\IN$:
  \[ f_n(y) = u_n(\xi(y))g(y).
  \]
  For the $\Cap_q$-polar set $N:=\kle{\xi^{-1}(P)\cap\klg{y\in\Omega_2:g(y)>0}}\cup N_1\cup N_2$ and $y\in\Omega_2\setminus N$ we get that 
  \[ f(y) = \lim_n f_n(y) = \lim_n u_n(\xi(y))g(y) = u(\xi(y))g(y).
  \]
\end{proof}

\subsection{Sobolev Spaces with Non-Vanishing Boundary Values: Local}\label{ss:h-nonzero}

Let $\Omega_1,\Omega_2\subset\IR^N$ be non-empty open sets and let $p,q\in (1,\infty)$. In this
subsection we assume that $T:\tsW^{1,p}(\Omega_1)\to\sW^{1,q}(\Omega_2)$ is a lattice homomorphism.
It follows from
\begin{itemize}
  \item Example \ref{ex:w1preg} that $\sW^{1,q}(\Omega_2)$ is regularizable,
  \item Lemma \ref{lem:h1ohp} that $T$ is $\overline\Omega_1$-local and positive (and continuous),
  \item Corollary \ref{cor:testrich1} that $\tsW^{1,p}(\Omega_1)\cap C_c(\overline\Omega_1)$ is rich in $\overline\Omega_1$ and
  \item Theorem \ref{thm:art} that $T|_{\tsW^{1,p}(\Omega_1)\cap C_c(\overline\Omega_1)}$ has a CoMu-Representation $(\xi,g)$,
\end{itemize}
that is, for all $u\in\tsW^{1,p}(\Omega_1)\cap C_c(\overline\Omega_1)$ and $f\in Tu$ there exists
a $\Cap_q$-polar set $N$ such that
\begin{equation} \label{eq:ident2}
  f(y) = u(\xi(y))g(y)\quad\mbox{ for all } y\in\Omega_2\setminus N.
\end{equation}
Note that here $\xi:\Omega_2\to\overline\Omega_1$. In order to extend Equation \eqref{eq:ident2} to
$u\in\su\in\tsW^{1,p}(\Omega_1)$ we need the following lemmata.

\begin{lemma}\label{lem:parta}
  Let $K_2\subset\Omega_2$ be a compact set and let $G_m:=\klg{y\in\Omega_2:g(y)>1/m}$ for $m\in\IN$.
  Then there exists a constant $C=C(K_2)$ such that for every compact set $K_1\subset\overline\Omega_1$ the
  following estimate holds:
  \[ \Cap_q(\xi^{-1}(K_1)\cap G_m\cap K_2) \leq
     C^qm^q\norm{T}^q\Cap_{p,\Omega_1}(K_1)^{q/p}.
  \]
\end{lemma}

\begin{proof}
  Let $\psi_{1,n}\in\tsW^{1,p}(\Omega_1)\cap C_c(\overline\Omega_1)$ and $\psi_2\in\cD(\Omega_2)$ be such
  that $\psi_{1,n}\geq 1$ on $K_1$, $\psi_2\geq 1$ on $K_2$ and
  $\norm{\psi_{1,n}}^p_{\sW^{1,p}(\Omega_1)}\leq \Cap_{p,\Omega_1}(K_1)+1/n$
  (see Theorem \ref{thm:biegert-cpt}).
  Let $f_n\in T\psi_{1,n}$ be fixed and let $N$ be a $\Cap_q$-polar set such that
  \[ f_n(y)=\psi_{1,n}(\xi(y))g(y)\quad\mbox{ for all }y\in\Omega_2\setminus N,n\in\IN.
  \]
  Then for $y\in\xi^{-1}(K_1)\cap G_m\cap K_2\cap N^c$ we get that
  \[ m\psi_2(y)f_n(y) \geq mg(y)\geq 1.
  \]
  Hence (using that $m\psi_2f_n$ is $\Cap_q$-quasi continuous) we get that
  \begin{eqnarray*}
    \Cap_q(\xi^{-1}(K)\cap G_m\cap K_2) &\leq& 
        \norm{mf_n\psi_2}^q_{\sW^{1,q}(\Omega_2)} \leq
         m^qC_2^q\norm{T}^q\norm{\psi_{1,n}}^q_{\sW^{1,p}(\Omega_1)} \\
        &\leq& m^qC_2^q\norm{T}^q\klr{\Cap_{p,\Omega_1}(K)+1/n}^{q/p}.
  \end{eqnarray*}
  For $n\to\infty$ the claim follows.
\end{proof}

\begin{lemma}\label{lem:partb}
  Let $K_2\subset\Omega_2$ be a compact set and let $G_m:=\klg{y\in\Omega_2:g(y)>1/m}$
  for $m\in\IN$. Then there exists a constant $C=C(K_2)$ such that for all sets $M$
  in $\overline\Omega_1$ the following estimate holds:
  \[ \Cap_q(\xi^{-1}(M)\cap G_m\cap K_2)\leq C^qm^q\norm{T}^q\Cap_{p,\Omega_1}(M)^{q/p}.
  \]
\end{lemma}

\begin{proof}
  Let $\varepsilon>0$. Then there exists an open set $O$ in the metric space
  $\overline\Omega_1$ containing $M$ such that $\Cap_{p,\Omega_1}(O)\leq\Cap_{p,\Omega_1}(M)+\varepsilon$
  (see Theorem \ref{thm:choquet}).
  Let $C_n\subset O$ be an increasing sequence of compact sets such that $\bigcup_n C_n=O$.
  By Lemma \ref{lem:parta}, using that $\Cap_{p,\Omega_1}$ is a Choquet Capacity
  (see Theorem \ref{thm:choquet}), we get that
  \begin{eqnarray*}
      \Cap_q(\xi^{-1}(M)\cap G_m\cap K_2) 
        &\leq& \Cap_q(\xi^{-1}(O)\cap G_m\cap K_2) \\
        &=& \lim_n \Cap_q(\xi^{-1}(C_n)\cap G_m\cap K_2) \\
        &\leq& \lim_n C^qm^q\norm{T}^q\Cap_{p,\Omega_1}(C_n)^{q/p}\\
        &=& C^qm^q\norm{T}^q\Cap_{p,\Omega_1}(O)^{q/p} \\
        &\leq& C^qm^q\norm{T}^q\klek{\Cap_{p,\Omega_1}(M)+\varepsilon}^{q/p}.
  \end{eqnarray*}
  For $\varepsilon\to 0+$ the claim follows.
\end{proof}

\begin{lemma}\label{lem:polarx}
  The set $\xi^{-1}(P)\cap\klg{y\in\Omega_2:g(y)>0}$ is $\Cap_q$-polar for every $\Cap_{p,\Omega_1}$-polar
  set $P\subset\overline\Omega_1$.
\end{lemma}

\begin{proof}
  Let $\omega_n\subset\Omega_2$ be an increasing sequence of compact sets such that
  $\bigcup_n \omega_n=\Omega_2$. It follows from Lemma \ref{lem:partb} that
  $\xi^{-1}(P)\cap G_m\cap \omega_n$ is $\Cap_q$-polar for all $m,n\in\IN$. The claim
  follows now from the identity
  \[ \xi^{-1}(P)\cap\klg{y\in\Omega_2:g(y)>0}=\bigcup_{n,m} \xi^{-1}(P)\cap G_m\cap\omega_n.
  \]
\end{proof}

\begin{theorem}\label{thm:nonzero}
  Let $\Omega_1,\Omega_2\subset\IR^N$ be non-empty open sets and let $p,q\in (1,\infty)$.
  Assume that $T:\tsW^{1,p}(\Omega_1)\to\sW^{1,q}(\Omega_2)$ is a lattice homomorphism.
  Then there exists a CoMu-Representation $(\xi,g)$ of $T$
  with $\xi:\Omega_2\to\overline\Omega_1$
  and $g:\Omega_2\to [0,\infty)$, that is,
  \[ T\su=(\su\circ\xi)g\quad\mbox{ for all }\su\in\tsW^{1,p}(\Omega_1).
  \]
  More precisely, this means that for every $u\in\su\in\tsW^{1,p}(\Omega_1)$ and every $f\in T\su$ there exists
  a $\Cap_q$-polar set $N\subset\Omega_2$ such that
  \[ f(y) = u(\xi(y))g(y)\quad\mbox{ for all } y\in\Omega_2\setminus N.
  \]
  Note that every $u\in\tsW^{1,p}(\Omega_1)$ has a unique trace on $\partial\Omega_1$ up to a
  $\Cap_{p,\Omega_1}$-polar set.
\end{theorem}

\begin{proof}
  Let $u_n\in\tsW^{1,p}(\Omega_1)\cap C_c(\overline\Omega_1)$ be a sequence of continuous functions
  converging in $\tsW^{1,p}(\Omega_1)$ to $\su$. For $u\in\su$ (after passing to a subsequence, Theorem \ref{thm:biegert-new})
  there exists a $\Cap_{p,\Omega_1}$-polar set $P$ such that $u_n\to u$ everywhere on $\overline\Omega_1\setminus P$.
  Now let $f_n\in Tu_n$ and $f\in T\su$ be fixed.
  Then (after passing to a subsequence, Theorem \ref{thm:biegert-new}) there exists a $\Cap_q$-polar set
  $N_1$ such that $f_n\to f$ everywhere on $\Omega_2\setminus N_1$. Let $N_2$ be a $\Cap_q$-polar
  set such that the following holds for all $y\in\Omega_2\setminus N_2$ and all $n\in\IN$:
  \[ f_n(y) = u_n(\xi(y))g(y).
  \]
  For the $\Cap_q$-polar set $N:=\kle{\xi^{-1}(P)\cap\klg{y\in\Omega_2:g(y)>0}}\cup N_1\cup N_2$ 
  (see Lemma \ref{lem:polarx}) and $y\in\Omega_2\setminus N$ we get that
  \[ f(y) = \lim_n f_n(y) = \lim_n u_n(\xi(y))g(y)=u(\xi(y))g(y).
  \]
\end{proof}

\subsection{Sobolev Spaces with Non-Vanishing Boundary Values: Global}\label{ss:h-nonzero2}

In this subsection we assume that $\Omega_1,\Omega_2\subset\IR^N$ are non-empty open sets,
$\Omega_1$ is bounded, $p,q\in (1,\infty)$ and $T:\tsW^{1,p}(\Omega_1)\to\tsW^{1,q}(\Omega_2)$
is a lattice homomorphism. Then there exists a CoMu-Representation $(\xi,g)$ of
$T:\tsW^{1,p}(\Omega_1)\to\sW^{1,q}(\Omega_2)$
(see Subsection \ref{ss:h-nonzero}) with $\xi:\Omega_2\to\overline\Omega_1$ and $g:\Omega_2\to[0,\infty)$,
that is, for all $u\in\tsW^{1,p}(\Omega_1)\cap C_c(\overline\Omega_1)$
\[ Tu = (u\circ\xi)g\qquad \mbox{$\Cap_{q,\Omega_2}$-quasi everywhere on }\Omega_2.
\]
Note that the $\Cap_q$- and $\Cap_{q,\Omega_2}$-polar sets in $\Omega_2$ coincide.

\begin{proposition}\label{prop:nonzero2}
  Under the above assumptions there exists a CoMu-Representation $(\xi^\star,g^\star)$ of $T$ with
  $\xi^\star:\overline\Omega_2\to\overline\Omega_1$ and $g^\star:\overline\Omega_2\to[0,\infty)$ such that
  for all $u\in\tsW^{1,p}(\Omega_1)\cap C_c(\overline\Omega_1)$
  \[ Tu = (u\circ \xi^\star)g^\star\quad\mbox{ $\Cap_{q,\Omega_2}$-quasi everywhere on $\overline\Omega_2$}.
  \]
\end{proposition}

\begin{proof}
  Since $g\in T1\in\tsW^{1,q}(\Omega_2)$ we get that $g:\Omega_2\to[0,\infty)$ has a unique extension
  $g^\star:\overline\Omega_2\to[0,\infty)$ which is $\Cap_{q,\Omega_2}$-quasi continuous on
  $\overline\Omega_2$. Since $\xi_jg=Tx_j\in\tsW^{1,q}(\Omega_2)$ we get that
  $\xi:\Omega_2\to\overline\Omega_1$ has an extension $\xi^\star:\overline\Omega_2\to\IR^N$
  which is $\Cap_{q,\Omega_2}$-quasi continuous on $\klg{y\in\overline\Omega_2:g^\star(y)\not=0}$.
  To see that $\xi^\star$ may be chosen such that $\xi^\star(\overline\Omega_2)\subset\overline\Omega_1$
  we let $u^\star\in C_b(\IR^N)$ be such that $u^\star>0$ on $\IR^N\setminus\overline\Omega_1$ and
  $u^\star=0$ on $\overline\Omega_1$. Then there exists a $\Cap_{q,\Omega_2}$-polar set $Q\subset\Omega_2$
  such that for $f:=0\in T(u^\star|_{\Omega_1})$
  \begin{equation}\label{eq:rep} 
     0 = u^\star(\xi^\star(y))g^\star(y)
  \end{equation}
  for all $y\in\Omega_2\setminus Q$.
  Since $f=0$ and $(u^\star\circ\xi^\star)g^\star$ are $\Cap_{q,\Omega_2}$-quasi continuous
  on $\overline\Omega_2$ we get by Theorem \ref{thm:biegert} (Uniqueness of the quasi continuous version)
  that there exists a $\Cap_{q,\Omega_2}$-polar set $Q^\star\subset\overline\Omega_2$ such that
  equation \eqref{eq:rep} holds for all $y\in\overline\Omega_2\setminus Q^\star$.
  Moreover, it follows that $(\xi^\star)^{-1}(\IR^N\setminus \overline\Omega_1)\cap\klg{y\in\overline\Omega_2:g(y)\not=0}\subset Q^\star$
  is a $\Cap_{q,\Omega_2}$-polar subset of $\overline\Omega_2$. Hence by changing $\xi^\star$ on a
  $\Cap_{q,\Omega_2}$-polar set and on $\klg{g=0}$ we get that
  $\xi^\star(\overline\Omega_2)\subset\overline\Omega_1$. Now let
  $u\in\tsW^{1,p}(\Omega_1)\cap C_c(\overline\Omega_1)$ and $f\in Tu$.
  Then there exists a $\Cap_{q,\Omega_2}$-polar set $Q\subset\Omega_2$ such that
  \begin{equation}\label{eq:rep2}
     f = (u\circ \xi^\star)g^\star
  \end{equation}
  everywhere on $\Omega_2\setminus Q$. Since both sides of equation \eqref{eq:rep2} are
  $\Cap_{q,\Omega_2}$-quasi continuous on $\overline\Omega_2$ this identity extends to hold
  $\Cap_{q,\Omega_2}$-quasi everywhere on $\overline\Omega_2$.
\end{proof}

In order to prove that the representation from Proposition \ref{prop:nonzero2} holds
even for all $\su\in\tsW^{1,p}(\Omega_1)$ we need the following lemmata.

\begin{lemma}\label{lem:part2a}
  For $m\in\IN$ let $G_m:=\klg{y\in\overline\Omega_2:g^\star(y)>1/m}$.
  Then for every compact set $K\subset\overline\Omega_1$ the following estimate holds:
  \[ \Cap_{q,\Omega_2}((\xi^\star)^{-1}(K)\cap G_m) \leq
     m^q\norm{T}^q\Cap_{p,\Omega_1}(K)^{q/p}.
  \]
\end{lemma}

\begin{proof}
  For $n\in\IN$ let $\psi_n\in\tsW^{1,p}(\Omega_1)\cap C_c(\overline\Omega_1)$ be such that
  $\psi_n\geq 1$ on $K$ and $\norm{\psi_n}^p_{\sW^{1,p}(\Omega_1)}\leq\Cap_{p,\Omega_1}(K)+1/n$.
  Let $f_n\in T\psi_n$
  be fixed and let $N\subset\overline\Omega_2$ be a $\Cap_{q,\Omega_2}$-polar set such that
  \[ f_n(y) = \psi_n(\xi^\star(y))g^\star(y)\qquad\mbox{ for all } y\in\overline\Omega_2\setminus N,n\in\IN.
  \]
  Then for $y\in(\xi^\star)^{-1}(K)\cap G_m\cap N^c$ we get that
  \[ mf_n(y) \geq mg^\star(y) \geq 1.
  \]
  Hence (using that $mf_n$ is $\Cap_{q,\Omega_2}$-quasi continuous) we get by Theorem
  \ref{thm:biegert-cap} that
  \begin{eqnarray*}
     \Cap_{q,\Omega_2}((\xi^\star)^{-1}(K)\cap G_m) &\leq&
      \norm{mf_n}^q_{\sW^{1,q}(\Omega_2)} \leq m^q\norm{T}^q\norm{\psi_{1,n}}^q_{\sW^{1,p}(\Omega_1)} \\
       &\leq& m^q\norm{T}^q\kle{\Cap_{p,\Omega_1}(K)+1/n}^{q/p}.
  \end{eqnarray*}
  For $n\to\infty$ the claim follows.
\end{proof}

\begin{lemma}\label{lem:part2b}
  For $m\in\IN$ let $G_m:=\klg{y\in\overline\Omega_2:g^\star(y)>1/m}$.
  Then for every set $M\subset\overline\Omega_1$ the following estimate holds:
  \[ \Cap_{q,\Omega_2}((\xi^\star)^{-1}(M)\cap G_m) \leq
     m^q\norm{T}^q\Cap_{p,\Omega_1}(M)^{q/p}.
  \]
\end{lemma}

\begin{proof}
  Let $\varepsilon>0$. Then there exists an open set $O$ in the metric space $\overline\Omega_1$
  containing $M$ such that $\Cap_{p,\Omega_1}(O)\leq\Cap_{p,\Omega_1}(M)+\varepsilon$.
  Let $C_n\subset O$ be an increasing sequence of compact sets such that $\bigcup_n C_n=O$.
  Now we get from Lemma \ref{lem:part2a} that
  \begin{eqnarray*}
    \Cap_{q,\Omega_2}((\xi^\star)^{-1}(M)\cap G_m) 
      &\leq& \Cap_{q,\Omega_2}\klr{(\xi^\star)^{-1}(O)\cap G_m} \\
      &=& \lim_n \Cap_{q,\Omega_2}\klr{(\xi^\star)^{-1}(C_n)\cap G_m} \\
      &\leq& m^q\norm{T}^q \lim_n \Cap_{p,\Omega_1}\klr{C_n}^{q/p} \\
      &=& m^q\norm{T}^q \Cap_{p,\Omega_1}\klr{O}^{q/p} \\
      &\leq& m^q\norm{T}^q \kle{\Cap_{p,\Omega_1}\klr{M}+\varepsilon}^{q/p}.
  \end{eqnarray*}
  For $\varepsilon\to 0+$ the claim follows.
\end{proof}

\begin{lemma}\label{lem:part2c}
  The set $(\xi^\star)^{-1}(P)\cap\klg{y\in\overline\Omega_2:g^\star(y)>0}$ is $\Cap_{q,\Omega_2}$-polar
  for every $\Cap_{p,\Omega_1}$-polar set $P\subset\overline\Omega_1$.
\end{lemma}

\begin{proof}
  For $m\in\IN$ let $G_m:=\klg{y\in\overline\Omega_2:g^\star(y)>1/m}$.
  Then the claim follows from Lemma \ref{lem:part2b} and the identity
  \[ (\xi^\star)^{-1}(P)\cap \klg{y\in\overline\Omega_2:g^\star(y)>0}=\bigcup_{m} (\xi^\star)^{-1}(P)\cap G_m.
  \]
\end{proof}

\begin{theorem}\label{thm:nonzero2}
  Let $\Omega_1,\Omega_2\subset\IR^N$ be non-empty open sets and let $p,q\in(1,\infty)$. Assume that
  $\Omega_1$ is bounded and that $T:\tsW^{1,p}(\Omega_1)\to\tsW^{1,q}(\Omega_2)$ is a lattice
  homomorphism. Then there exists a CoMu-Representation $(\xi^\star,g^\star)$ of $T$
  with $\xi^\star:\overline\Omega_2\to\overline\Omega_1$ and $g:\overline\Omega_2\to[0,\infty)$ such that
  for all $\su\in\tsW^{1,p}(\Omega_1)$
  \[ T\su = (\su\circ\xi^\star)g^\star\qquad\mbox{ $\Cap_{q,\Omega_2}$-quasi everywhere on $\overline\Omega_2$}.
  \]
  More precisely, this means that for every $u\in\su\in\tsW^{1,p}(\Omega_1)$ and every $f\in T\su$ there exists
  a $\Cap_{q,\Omega_2}$-polar set $N\subset\overline\Omega_2$ such that
  \[ f(y) = u(\xi(y))g(y)\quad\mbox{ for all } y\in\overline\Omega_2\setminus N.
  \]
\end{theorem}

\begin{proof}
  Let $u\in\su\in\tsW^{1,p}(\Omega_1)$. Then there exist $u_n\in\tsW^{1,p}(\Omega_1)\cap C_c(\overline\Omega_1)$ 
  and a $\Cap_{p,\Omega_1}$-polar set $P$ such that $u_n\to\su$ in $\tsW^{1,p}(\Omega_1)$ and
  $u_n\to u$ everywhere on $\overline\Omega_1\setminus P$. Now let $f_n\in Tu_n$ and $f\in T\su$.
  Then (after passing to a subsequence) there exists a $\Cap_{q,\Omega_2}$-polar set $N_1$ such that
  $f_n\to f$ everywhere on $\overline\Omega_2\setminus N_1$. Let $N_2$ be a
  $\Cap_{q,\Omega_2}$-polar set such that the following holds for all
  $y\in\overline\Omega_2\setminus N_2$ and all $n\in\IN$
  \[ f_n(y) = u_n(\xi^\star(y))g^\star(y).
  \]
  For the $\Cap_{q,\Omega_2}$-polar set $N:=\kle{(\xi^\star)^{-1}(P)\cap\klg{y\in\overline\Omega_2:g(y)>0}}\cup N_1\cup N_2$
  (see Lemma \ref{lem:part2c}) and $y\in\overline\Omega_2\setminus N$ we get that
  \[ f(y) = \lim_n f_n(y) = \lim_n u_n(\xi^\star(y))g^\star(y) = u(\xi^\star(y))g^\star(y).
  \]
\end{proof}

% Now let $u\in\su\in \tsW^{1,p}(\Omega_1)$
%  be fixed. Then there exists $u_n\in\tsW^{1,p}(\Omega_1)\cap C_c(\overline\Omega_1)$ and a
%  $\Cap_{p,\Omega_1}$-polar set $S\subset\overline\Omega_1$ such that
%  $u_n\to u$ in $\tsW^{1,p}(\Omega_1)$ and everywhere on $\overline\Omega_1\setminus S$.
%  Let $f_n\in Tu_n$ and $f\in T\su$ be fixed. Then (after subsequence) there exists a
%  $\Cap_{q,\Omega_2}$-polar set $R$ such that $f_n\to f$ everywhere on $\overline\Omega_2\setminus R$.
%  Denote by $R_n\subset\overline\Omega_2$ a $\Cap_{q,\Omega_2}$-polar set such that
%  equation \eqref{eq:rep2} holds for $f_n$ and $u_n$ everywhere on $\overline\Omega_2\setminus R_n$.
%  Then for the $\Cap_{q,\Omega_2}$-polar set $R:=(\xi^\star)^{-1}(S)\cup \klr{\bigcup_n R_n}$
%\end{proof}

\subsection{Sobolev Spaces with Vanishing Boundary Values: Lattice Isomorphisms}

In this subsection we assume that $p,q\in (1,\infty)$, $\Omega_1,\Omega_2\subset\IR^N$
are non-empty open sets and $T:\sW^{1,p}_0(\Omega_1)\to\sW^{1,q}(\Omega_2)$ is a lattice homomorphism.
Let $(\xi,g)$ be a CoMu-Representation of $T$ which exists by Theorem \ref{thm:zero}.

\begin{definition}
  For $p\in(1,\infty)$ and $N\subset\Omega$ we define the Banach space $\sW^{1,p}_0(\Omega,N)$ by
  \[ \sW^{1,p}_0(\Omega,N) := \klg{\su\in\sW^{1,p}(\Omega):\su=0\;\Cap_p\mbox{-q.e. on }N},\quad
     \norm{\su}_{\sW^{1,p}(\Omega,N)}:=\norm{\su}_{\sW^{1,p}(\Omega)}.
  \]
\end{definition}

\begin{proposition}\label{prop:zero}
  Let $p\in(1,\infty)$ and let $N\subset\Omega$ be an arbitrary set. Then
  \[ \cD(\Omega) \subset \sW^{1,p}_0(\Omega,N)\quad\mbox{ if and only if }\quad
     \Cap_p(N)=0.
  \]
\end{proposition}

\begin{proof}
  Let $\omega_k\subset\subset\Omega$ be open sets such that $\bigcup_k\omega_k=\Omega$
  and let $\varphi_k\in\cD(\Omega)$ be such that $\varphi_k\equiv 1$ on $\omega_k$.
  Assume now that $\cD(\Omega)\subset\sW^{1,p}_0(\Omega,N)$.
  Then $\varphi_k\in\sW^{1,p}_0(\Omega,N)$ and hence $\Cap_p(N\cap\omega_k)=0$.
  This shows that $\Cap_p(N)=0$. Assume now that $\Cap_p(N)=0$. Then
  $\sW^{1,p}_0(\Omega,N)=\sW^{1,p}(\Omega)$ and hence $\cD(\Omega)\subset\sW^{1,p}_0(\Omega,N)$.
\end{proof}

\begin{proposition}\label{prop:postest}
  If $\cD(\Omega_2)\subset T\sW^{1,p}_0(\Omega_1)$ then $g$ is strictly positive
  $\Cap_q$-q.e. on $\Omega_2$.
\end{proposition}

\begin{proof}
  Let $N:=\klg{y\in\Omega_2:g(y)=0}$. Then $\cD(\Omega_2)\subset T\sW^{1,p}_0(\Omega_1)\subset \sW^{1,q}_0(\Omega_2,N)$.
  Therefore $\Cap_q(N)=0$ by Proposition \ref{prop:zero}, that is,
  $g>0$ $\Cap_q$-q.e. on $\Omega_2$.
\end{proof}

\begin{definition}
  For $j=1,2$ let $\topo_j$ be a topological space and $\Cap_j$ be a Choquet capacity on $\topo_j$.
  Then a mapping $\tau:\topo_1\to\topo_2$ is called {\em ${\Cap_1}$-${\Cap_2}$-quasi invertible} 
  if there exist a $\Cap_1$-polar set $S$ and a $\Cap_2$-polar set $R$ such that
  $\tau:\topo_1\setminus S\to \topo_2\setminus R$ is bijective. In this case we let
  $\tau^{-1}:\topo_2\to\topo_1$ be given by
  \[ \tau^{-1}(y) := 
      \begin{cases}
       \klrk{\tau|_{\topo_1\setminus S}}^{-1}(y) \mbox{ if }y\in \topo_2\setminus R\\
       \in \topo_2\mbox{ arbitrarily if } y\in R.
      \end{cases}
  \]
\end{definition}

\begin{theorem}
  Let $\Omega_1,\Omega_2\subset\IR^N$ be non-empty open sets and let $p,q\in (1,\infty)$.
  Assume that $T:\sW^{1,p}_0(\Omega_1)\to \sW^{1,q}_0(\Omega_2)$ is a lattice isomorphism.
  Then there exists a CoMu-Representation $(\xi,g)$ of $T$ with $\xi:\Omega_2\to\Omega_1$
  and $g:\Omega_2\to (0,\infty)$ such that $\xi$ is $\Cap_q$-$\Cap_p$-quasi invertible and
  $(\xi^{-1},1/g\circ \xi^{-1})$ is a CoMu-Representation for $T^{-1}$.
\end{theorem}

\begin{remark}
  In the above theorem it does not matter whether such lattice isomorphisms exist
  (for $p\not=q$) or not.
\end{remark}

\begin{proof}
  Note that the inverse $T^{-1}$ of $T$ is again a lattice homomorphism
  (see Aliprantis and Burkinshaw \cite[Theorem 7.3]{aliprantis:85:po}).
  By Theorem \ref{thm:zero} and Proposition \ref{prop:postest} there
  exist CoMu-Representations $(\xi,g)$ of $T$ and $(\eta,h)$ of  $T^{-1}$ with
  $g(\Omega_2)\subset(0,\infty)$ and $h(\Omega_1)\subset(0,\infty)$:
  \[ T\su=(\su\circ\xi)g\;\mbox{ for all } \su\in\sW^{1,p}_0(\Omega_1)\quad\mbox{ and }\quad
     T^{-1}\sv=(\sv\circ\eta)h\;\mbox{ for all } \sv\in\sW^{1,q}_0(\Omega_2).
  \]
  Let $v\in\sv\in\sW^{1,q}_0(\Omega_2)$, $u\in\su:=T^{-1}\sv$ and $w\in T\su=\sv$.
  Then there exist a $\Cap_p$-polar set $P_1$ and a $\Cap_q$-polar set $Q_1$ such that
  $u(x)=v(\eta(x))h(x)$ for all $x\in\Omega_1\setminus P_1$ and
  $v(y) = w(y) = u(\xi(y))g(y)$ for all $y\in\Omega_2\setminus Q_1$.
  Hence we conclude that
  \[ v(y) = v(\eta(\xi(y)))h(\xi(y))g(y)\quad\mbox{ for all } y\in\Omega_2\setminus(Q_1\cup\xi^{-1}(P_1)).
  \]
  It follows from Lemma \ref{lem:polar} and Proposition \ref{prop:postest} that the set
  $Q:=Q_1\cup\xi^{-1}(P_1)$ is a $\Cap_q$-polar set.
  Now let $\omega_n\subset\subset\Omega_2$ be a sequence of open sets such that $\bigcup_n \omega_n=\Omega_2$
  and let $v_n,w_{n,j}\in\cD(\Omega_2)$ be such that $v_n(y)=1$ and $w_{n,j}(y)=y_j$ for all $y=(y_1,\dots,y_N)^t\in\omega_n$.
  Then there exists a $\Cap_q$-polar set $R$ such that for all $y\in\Omega_2':=\Omega_2\setminus R$, all $n\in\IN$
  and all $j\in\klg{1,\dots,N}$
  \[ v_n(y) = v_n(\eta(\xi(y)))h(\xi(y))g(y)\;\quad\mbox{ and }\quad
     w_{n,j}(y) = w_{n,j}(\eta(\xi(y)))h(\xi(y))g(y).
  \]
  Let $y\in\Omega_2'$ be fixed and let $n\in\IN$ be such that $y\in\omega_n$ and $\eta(\xi(y))\in\omega_n$. Then we get that
  \begin{eqnarray*}
     1 &=& v_n(y) = v_n(\eta(\xi(y)))h(\xi(y))g(y) = h(\xi(y))g(y);\\
     y_j &=& w_n(y) = w_n(\eta(\xi(y)))h(\xi(y))g(y) = w_n(\eta(\xi(y))) = \eta(\xi(y))_j.
  \end{eqnarray*}
  Therefore $g=1/(h\circ\xi)$ everywhere on $\Omega_2'$, $\xi:\Omega_2'\to\xi(\Omega_2')$ is bijective
  and $\eta:\xi(\Omega_2')\to \Omega_2'$ is its inverse.
  Interchanging the role of $T$ and $T^{-1}$ we get a $\Cap_p$-polar set $S$ and a set
  $\Omega_1':=\Omega_1\setminus S$ such that
  $\eta:\Omega_1'\to\eta(\Omega_1')$ is bijective and $\xi:\eta(\Omega_1')\to\Omega_1'$ is
  its inverse. For $\tOm_1:=\Omega_1'\cup\xi(\Omega_2')$ and $\tOm_2:=\eta(\Omega_1')\cup\Omega_2'$
  we get that $\xi:\tOm_2\to\tOm_1$ is bijective and $\eta:\tOm_1\to\tOm_2$ is its inverse.
\end{proof}

\subsection{Sobolev Spaces with Non-Vanishing Boundary Values: Lattice Isomorphisms}

In this subsection we assume that $p,q\in (1,\infty)$, $\Omega_1,\Omega_2\subset\IR^N$
are bounded non-empty open sets and $T:\tsW^{1,p}(\Omega_1)\to\tsW^{1,q}(\Omega_2)$
is a lattice isomorphism.
%Whether such lattice isomorphisms exist (for $p\not=q$) or not will not be considered here.
Let $(\xi^\star,g^\star)$ be a CoMu-Representation of $T$ which exists by Theorem \ref{thm:nonzero2}.

\begin{lemma}\label{lem:pos2}
  The function $g^\star$ is strictly positive $\Cap_{q,\Omega_2}$-quasi everywhere on
  $\overline\Omega_2$.
\end{lemma}

\begin{proof}
  Let $N:=\klgk{y\in\overline\Omega_2:g^\star(y)=0}$ and let $\omega_n\subset\overline\Omega_2$ be
  a sequence of compact sets such that $\bigcup_n \omega_n=\overline\Omega_2$.
  Let $\varphi_n\in\cD(\IR^N)$ be such that $\varphi_n\geq 1$ on $\omega_n$.
  Since $\varphi_n|_{\overline\Omega_2}\in\tsW^{1,q}(\Omega_2)$ is in the image of
  $T$, we get that
  $\Cap_{q,\Omega_2}(\omega_n\cap N)=0$. Therefore $\Cap_{q,\Omega_2}(N)=0$.
\end{proof}

\begin{theorem}
  Let $\Omega_1,\Omega_2\subset\IR^N$ be bounded non-empty open sets and let
  $p,q\in (1,\infty)$. Assume that $T:\tsW^{1,p}(\Omega_1)\to \tsW^{1,q}(\Omega_2)$ is a lattice
  isomorphism. Then there exists a CoMu-Representation $(\xi^\star,g^\star)$ of $T$ with
  $\xi^\star:\overline\Omega_2\to\overline\Omega_1$ and $g^\star:\overline\Omega_2\to (0,\infty)$
  such that $\xi^\star$ is $\Cap_{q,\Omega_2}$-$\Cap_{p,\Omega_1}$-quasi invertible and
  $\klr{(\xi^\star)^{-1},1/g^\star\circ (\xi^\star)^{-1}}$ is a CoMu-Representation for $T^{-1}$.
\end{theorem}

\begin{proof}
  Note that the inverse $T^{-1}$ of $T$ is again a lattice homomorphism
  (see Aliprantis and Burkinshaw \cite[Theorem 7.3]{aliprantis:85:po}).
  By Theorem \ref{thm:nonzero2} and Lemma \ref{lem:pos2} there
  exist CoMu-Representations $(\xi^\star,g^\star)$ of $T$ and $(\eta^\star,h^\star)$ of  $T^{-1}$ with
  $g^\star(\overline\Omega_2)\subset(0,\infty)$ and $h^\star(\overline\Omega_1)\subset(0,\infty)$, that is,
  \[ T\su=(\su\circ\xi^\star)g^\star\mbox{ for all } \su\in\tsW^{1,p}(\Omega_1)\quad\mbox{ and }\quad
     T^{-1}\sv=(\sv\circ\eta^\star)h^\star\mbox{ for all } \sv\in\tsW^{1,q}(\Omega_2).
  \]
  Let $v\in\sv\in\tsW^{1,q}(\Omega_2)$, $u\in\su:=T^{-1}\sv$ and $w\in T\su=\sv$.
  Then there exist a $\Cap_{p,\Omega_1}$-polar set $P_1$ and a $\Cap_{q,\Omega_2}$-polar set $Q_1$ such that
  $u(x)=v(\eta^\star(x))h^\star(x)$ for all $x\in\overline\Omega_1\setminus P_1$ and
  $v(y) = w(y) = u(\xi^\star(y))g^\star(y)$ for all $y\in\overline\Omega_2\setminus Q_1$.
  Hence we conclude that
  \[ v(y) = v(\eta^\star(\xi^\star(y)))h^\star(\xi^\star(y))g^\star(y)\quad\mbox{ for all }
    y\in\overline\Omega_2\setminus(Q_1\cup(\xi^\star)^{-1}(P_1)).
  \]
  It follows from Lemma \ref{lem:part2c} and Lemma \ref{lem:pos2} that the set
  $Q:=Q_1\cup(\xi^\star)^{-1}(P_1)$ is a $\Cap_{q,\Omega_2}$-polar set.
  Now let $v,w_j\in\tsW^{1,q}(\Omega_2)\cap C_c(\overline\Omega_2)$ be given by
  $v(y):=1$ and $w_j(y):=y_j$. Then there exists a $\Cap_{q,\Omega_2}$-polar set $R$ such that
  for all $y\in\Omega_2':=\overline\Omega_2\setminus R$ and all $j\in\klg{1,\dots,N}$
  \begin{eqnarray*}
     1&=&v(y) = v(\eta^\star(\xi^\star(y)))h^\star(\xi^\star(y))g^\star(y) = h^\star(\xi^\star(y))g^\star(y)\\
     y_j&=&w_j(y) = w_j(\eta^\star(\xi^\star(y)))h^\star(\xi^\star(y))g^\star(y) = w_j(\eta^\star(\xi^\star(y))).
  \end{eqnarray*}
  Therefore $g^\star=1/(h^\star\circ\xi^\star)$ everywhere on $\Omega_2'$,
  $\xi^\star:\Omega_2'\to\xi^\star(\Omega_2')$ is bijective and $\eta^\star:\xi^\star(\Omega_2')\to \Omega_2'$ is its inverse.
  Interchanging the role of $T$ and $T^{-1}$ we get a $\Cap_{p,\Omega_1}$-polar set $S$ and a set
  $\Omega_1':=\overline\Omega_1\setminus S$ such that
  $\eta^\star:\Omega_1'\to\eta(\Omega_1')$ is bijective and $\xi^\star:\eta^\star(\Omega_1')\to\Omega_1'$ is
  its inverse. For $\tOm_1:=\Omega_1'\cup\xi^\star(\Omega_2')$ and $\tOm_2:=\eta^\star(\Omega_1')\cup\Omega_2'$
  we get that $\xi^\star:\tOm_2\to\tOm_1$ is bijective and $\eta^\star:\tOm_1\to\tOm_2$ is its inverse.
\end{proof}

\subsection{$L^p$ Spaces}
Let $\Omega_1,\Omega_2\subset\IR^N$ be non-empty open sets and let $p,q\in[1,\infty]$.
In this section we assume that $T:L^p(\Omega_1)\to L^q(\Omega_2)$ is a lattice homomorphism.
It follows from
\begin{itemize}
  \item Example \ref{ex:lpreg} that $L^q(\Omega_2)$ is regularizable,
  \item Lemma \ref{lem:l2oh} that $T$ is $\Omega_1$-local and positive (and continuous),
  \item Corollary \ref{cor:testrich0} that $C_c(\Omega_1)$ is rich in $\Omega_1$ and
  \item Theorem \ref{thm:art} that $T|_{C_c(\Omega_1)}$ has a
        CoMu-Representation $(\xi,g)$,
\end{itemize} 
that is, for all $u\in C_c(\Omega_1)\supset\cD(\Omega_1)$ and $f\in Tu$ there exists
a Lebesgue nullset $N$ such that
\begin{equation}\label{eq:lp}
   f(y)=u(\xi(y))g(y)\quad\text{ for all }y\in\Omega_2\setminus N.
\end{equation}
In order to extend Equation \eqref{eq:lp} to $u\in L^p(\Omega_1)$ we need the following lemma.

\begin{lemma}
  The set $\xi^{-1}(P)\cap\klg{y\in\Omega_2:g(y)>0}$ is a Lebesgue nullset for every
  Lebesgue nullset $P\subset\Omega_1$.
\end{lemma}

\begin{proof}
  Using that for a compact set $K\subset\Omega_1$ we have that
  \[ \lambda(K)=\inf\klg{\norm{u}_{L^p(\Omega_1)}^p:u\in C_c(\Omega_1),u\geq 1\mbox{ on }K}
  \]
  we get that there are $\varphi_n\in C_c(\Omega_1)$ such that
  $\varphi_n\geq 1$ on $K$ and $\norm{\varphi_n}_{L^p(\Omega_1)}^p\leq\lambda(K)+1/n$.
  Let $G_m:=\klg{y\in\Omega_2:g(y)>1/m}$ and $f_n\in T\varphi_n$. Then
  there exists a Lebesgue nullset $P_2\subset\Omega_2$ such that
  \[ f_n(y) = \varphi_n(\xi(y))g(y)\geq 1/m\quad\mbox{ for all } y\in\xi^{-1}(K)\cap G_m\cap P_2^c.
  \]
  This shows (using the outer Lebesgue measure $\lambda^\star$) that 
  \[ \lambda^\star(\xi^{-1}(K)\cap G_m) \leq m^q\norm{f_n}_{L^q(\Omega)}^q\leq
     m^q\norm{T}^q\norm{\varphi_n}_{L^p(\Omega_1)}^q \leq
     m^q\norm{T}^q(\lambda(K)+1/n)^{q/p}.
  \]
  For $n\to\infty$ we get that
  \[ \lambda^\star(\xi^{-1}(K)\cap G_m)\leq m^q\norm{T}^q\lambda(K)^{q/p}.
  \]
  Now let $O\subset\Omega_1$ be an open and fixed set and let $K_j\subset O$ be an
  increasing sequence of compact sets such that $\bigcup_j K_j=O$. Then
  \[ \lambda^\star(\xi^{-1}(O)\cap G_m) =
     \lim_j \lambda^\star(\xi^{-1}(K_j)\cap G_m)\leq
     \lim_j m^q\norm{T}^q\lambda(K_j)^{q/p} = m^q\norm{T}^q\lambda(O)^{q/p}.
  \]
  Now take the Lebesgue nullset $P\subset\Omega_1$. Then for $\varepsilon>0$ there exists
  an open set $O$ containing $P$ such that $\lambda(O)\leq\varepsilon$. Hence
  \[ \lambda^\star(\xi^{-1}(P)\cap G_m) \leq
     \lambda^\star(\xi^{-1}(O)\cap G_m) \leq m^q\norm{T}^q\varepsilon^{q/p}.
  \]
  For $\varepsilon\to 0$ we get that $\xi^{-1}(P)\cap G_m$ is a Lebesgue nullset.
  Let $G:=\bigcup G_m$. Then the claim follows form the equality
  $\xi^{-1}(P)\cap G=\bigcup \xi^{-1}(P)\cap G_m$.
\end{proof}

\begin{theorem}\label{thm:mainlp}
  Let $\Omega_1,\Omega_2\subset\IR^N$ be non-empty open sets and let $p,q\in[1,\infty]$.
  Assume that $T:L^p(\Omega_1)\to L^q(\Omega_2)$ is a lattice homomorphism. Then there
  exist a CoMu-Representation $(\xi,g)$ of $T$ with $\xi:\Omega_1\to\Omega_2$ and
  $g:\Omega_2\to[0,\infty)$, that is,
  \[ T\su=(\su\circ\xi)g\quad\text{ for all }\su\in L^p(\Omega_1).
  \]
  More precisely, this means that for every $u\in\su\in L^p(\Omega_1)$ and every
  $f\in T\su$ there exists a Lebesgue nullset $N\subset\Omega_2$ such that
  \[ f(y) = u(\xi(y))g(y)\quad\mbox{ for all } y\in\Omega_2\setminus N.
  \]
\end{theorem}

\begin{proof}
  Let $u_n\in\cD(\Omega_1)$ be a sequence of test functions converging in $L^p(\Omega_1)$
  to $\su\in L^p(\Omega_1)$. For $u\in\su$ (after passing to a subsequence) there exists a
  Lebesgue nullset $P$ such that $u_n\to u$ everywhere on $\Omega_1\setminus P$.
  Now let $f_n\in Tu_n$ and $f\in T\su$ be fixed. Then (after passing to a subsequence)
  there exists a Lebesgue nullset $N_1$ such that $f_n\to f$ everywhere on
  $\Omega_2\setminus N_1$. Let $N_2$ be a Lebesgue nullset such that the following
  holds for all $y\in\Omega_2\setminus N_2$ and all $n\in\IN$:
  \[ f_n(y) = u_n(\xi(y))g(y).
  \]
  For the Lebesgue nullset $N:=\kle{\xi^{-1}(P)\cap\klg{y\in\Omega_2:g(y)>0}}\cup N_1\cup N_2$
  and $y\in\Omega_2\setminus N$ we get that
  \[ f(y) = \lim_n f_n(y) = \lim_n u_n(\xi(y))g(y) = u(\xi(y))g(y).
  \]
\end{proof}

Here we should mention the following representation theorem
of Feldman and Porter \cite[Theorem 1]{feldman:86:obl} for
lattice homomorphisms between certain Banach lattices.

\begin{theorem}(Feldman and Porter).
  Let $E$ and $F$ be Banach lattices having locally compact representation
  spaces $X$ and $Y$ respectively
  (see Schaefer \cite[Definition III.5.4]{schaefer:74:blp}) and let $T:E\to F$ be a lattice homomorphism
  satisfying $T(C_c(X))\subset C_c(Y)$. Then there are a non-negative
  function $g:Y\to\IR$ and a function $\xi:Y\to X$, both continuous on
  $P:=\klg{y\in Y:g(y)>0}$, such that for $u\in E$
  \[ (Tu)(y) = \begin{cases}
                 u(\xi(y))g(y) & \text{ if } y\in P \\
                 0 & \text{ if } y\in Y\setminus P.
               \end{cases}
  \]
\end{theorem}

\subsection{$L^p$ and Sobolev spaces}

\begin{theorem}
  Let $\Omega_1,\Omega_2\subset\IR^N$ be non-empty open sets where $\Omega_1$
  has finite Lebesgue measure and let $p\in[1,\infty)$ and $q\in[1,\infty]$.
  Assume that $T:W^{1,p}(\Omega_1)\to L^q(\Omega_2)$ is a lattice homomorphism.
  Then there exist a CoMu-Representation $(\xi,g)$ of $T$ with
  $\xi:\Omega_2\to\Omega_1$ and $g:\Omega_2\to[0,\infty)$, that is,
  \[ T\su = (\su\circ\xi)g\quad\text{ for all }\su\in W^{1,p}(\Omega_1).
  \]
  More precisely, this means that for every $u\in\su\in W^{1,p}(\Omega)$ and
  every $f\in T\su$ there exists a Lebesgue nullset $N\subset\Omega_2$ such that
  \[ f(y) = u(\xi(y))g(y)\quad\text{ for all } y\in\Omega_2\setminus N.
  \]
\end{theorem}

\begin{proof}
  Let $S$ be the restriction of $T$ to the sublattice $L:=W^{1,p}(\Omega_1)\cap L^\infty(\Omega_1)$.
  Then $S:L\to L^q(\Omega_2)$ is a lattice homomorphism and $L$ dominates
  $L^\infty(\Omega_1)$. Using that $L^q(\Omega_2)$ is complete
  vector lattice (see Meyer-Nieberg \cite[Example v, page 9]{meyer-nieberg:91:bl})
  we can extend $S$ (see Bernau \cite[Theorem 3.1]{bernau:86:evl}) to a vector lattice
  homomorphism $\tilde S:L^\infty(\Omega_1)\to L^q(\Omega_2)$.
  By Theorem \ref{thm:mainlp} we get that there exist a CoMu-Representation $(\xi,g)$
  of $\tilde S$, that is, $\tilde S\su=(\su\circ\xi)g$ for all
  $\su\in L^\infty(\Omega_1)$, in particular $T\su=S\su=\tilde S\su=(\su\circ\xi)g$
  for all $\su\in L=W^{1,p}(\Omega_1)\cap L^\infty(\Omega_1)$. For $\su\in W^{1,p}(\Omega_1)$
  define $\su_n:=(\su\vee (-n))\wedge n$. Then $\su_n\to\su$ almost everywhere on $\Omega_1$ and
  in $W^{1,p}(\Omega_1)$. Arguing as in the proof of Theorem \ref{thm:mainlp}
  the claim follows.
\end{proof}

\section{Examples and Remarks}\label{sec:examples}

\begin{remark}
  The assumption that $T$ is a lattice homomorphism in the previous theorems
  can be reduced to the condition that $T$ is a linear, order bounded and
  disjointness preserving operator. This follows from the following theorem
  (see \cite[Theorem 1.2]{bernau:92:sevl}).
  \begin{quote} Let $E$ and $F$ be vector lattices with $F$ Archimedean and
                $T$ an order bounded linear operator from $E$ into $F$
                such that $\betrag{T\su}\wedge\betrag{T\sv}=0$ for all
                $\su,\sv\in E$ with $\betrag{\su}\wedge\betrag{\sv}=0$. Then
                there exist lattice homomorphisms $T^+$, $T^-$ and $\betrag{T}$
                from $E$ into $F$ such that $T=T^+-T^-$, $(T^+)x=(Tx)^+$ and
                $(T^-)x=(Tx)^-$ $(0\leq x\in E)$, $\betrag{T}=T^++T^-$ and
                $\betrag{Tx}=\betrag{T}(\betrag{x})$ for $x\in E$.
  \end{quote}
\end{remark}

\begin{example}{\bf (Restriction to a smaller set)}
  Let $\Omega_1,\Omega_2\subset\IR^2$ be given by $\Omega_1:=(0,2)\times (0,2)$ and $\Omega_2:=(0,1)\times(0,1)$.
  Then $T:\sW^{1,p}(\Omega_1)\to\sW^{1,q}(\Omega_2)$ given by $Tu:=u|_{\Omega_2}$ is a surjective lattice homomorphism.
\end{example}

\begin{example}{\bf (Extension by zero)}
  Let $\Omega_1,\Omega_2\subset\IR^2$ be given by $\Omega_1:=(0,2)\times (0,1)$ and $\Omega_2:=(0,2)\times(0,2)$.
  Then $T:\sW^{1,p}_0(\Omega_1)\to\sW^{1,p}_0(\Omega_2)$ given by
  $Tu(x):=u(x)$ if $x\in\Omega_1$ and $Tu(x):=0$ if $x\in\Omega_2\setminus\Omega_1$ defines an injective
  lattice homomorphism.
\end{example}

\begin{example}{\bf (Reflection at the boundary)}
  Let $\Omega_1,\Omega_2\subset\IR^2$ be given by $\Omega_1:=(0,2)\times (0,1)$ and
  $\Omega_2:=(0,2)\times(0,2)$.
  Then $T:\sW^{1,p}(\Omega_1)=\tsW^{1,p}(\Omega_1)\to\sW^{1,p}(\Omega_2)$ given by
  $Tu(x):=u(x)$ if $x\in\Omega_1$ and $Tu(x):=u(x_1,2-x_2)$ if $x\in\Omega_2\setminus\overline\Omega_1$
  defines an injective
  lattice homomorphism. Let $\xi:\Omega_2\to\overline\Omega_1$ and $g:\Omega_2\to[0,\infty)$ be given from
  Theorem \ref{thm:nonzero}. Let $u(x):=x_2$. Then $0<u<1$ on $\Omega_1$ but
  $Tu=1$ $\Cap_p$-quasi everywhere on $(0,2)\times\klg{1}$.
  Using that $g=1$ $\Cap_p$-quasi everyhwere on $\Omega_2$ we get that the mapping
  $\xi:\Omega_2\to\overline\Omega_1$ from Theorem \ref{thm:nonzero} can not be changed
  on a $\Cap_{q,\Omega_2}$-polar set such that $\xi(\Omega_2)\subset\Omega_1$.
  This example shows also that the trace is needed to get a representation
  $\Cap_{p,\Omega_2}$-quasi every.
\end{example}

\bibliography{biblio2}

\providecommand{\mathbb}[1]{\mathbf{#1}}\providecommand{\cprime}{$'$}
\providecommand{\bysame}{\leavevmode\hbox to3em{\hrulefill}\thinspace}
\providecommand{\MR}{\relax\ifhmode\unskip\space\fi MR }
\providecommand{\MRhref}[2]{%
  \href{http://www.ams.org/mathscinet-getitem?mr=#1}{#2}
}
\providecommand{\href}[2]{#2}
\begin{thebibliography}{10}

\bibitem{adams:96:fsp}
David~R. Adams and Lars~Inge Hedberg, \emph{Function spaces and potential
  theory}, Grundlehren der mathematischen Wissenschaften, vol. 314,
  Springer-Verlag, Berlin, 1996. \MR{97j:46024}

\bibitem{roeckner:97:pof}
S.~Albeverio, Zhiming Ma, and Michael R{\"o}ckner, \emph{Partitions of unity in
  sobolev spaces over infinite dimensional state spaces.}, J. Funct. Anal.
  \textbf{143} (1997), no.~1, 247--268.

\bibitem{aliprantis:85:po}
Charalambos~D. Aliprantis and Owen Burkinshaw, \emph{Positive operators},
  Springer, Dordrecht, 2006, Reprint of the 1985 original. \MR{MR2262133}

\bibitem{arendt:84:rpo}
W.~Arendt, \emph{Resolvent positive operators and integrated semigroups},
  Semesterbericht Funktionalanalysis, T\"ubingen, Sommersemester 1984 (1984),
  73--101.

\bibitem{arendt:03:lrb}
Wolfgang Arendt and Mahamadi Warma, \emph{The {L}aplacian with {R}obin boundary
  conditions on arbitrary domains}, Potential Anal. \textbf{19} (2003), no.~4,
  341--363. \MR{1 988 110}

\bibitem{banach:32:tol}
Stefan Banach, \emph{Th\'eorie des op\'erations lin\'eaires}, \'Editions
  Jacques Gabay, Sceaux, 1993, Reprint of the 1932 original. \MR{MR1357166
  (97d:01035)}

\bibitem{bernau:86:evl}
S.~J. Bernau, \emph{Extension of vector lattice homomorphisms}, J. London Math.
  Soc. (2) \textbf{33} (1986), no.~3, 516--524. \MR{MR850967 (87h:47088)}

\bibitem{bernau:92:sevl}
\bysame, \emph{Sums and extensions of vector lattice homomorphisms}, Acta Appl.
  Math. \textbf{27} (1992), no.~1-2, 33--45, Positive operators and semigroups
  on Banach lattices (Cura\c cao, 1990). \MR{MR1184875 (93f:47041)}

\bibitem{biegert:08:trc}
Markus Biegert, \emph{The relative capacity}, Preprint, June 2008, available
  online at \url{http://arxiv.org/abs/0806.1417v3}.

\bibitem{bouleau:91:df}
Nicolas Bouleau and Francis Hirsch, \emph{Dirichlet forms and analysis on
  {W}iener space}, de Gruyter Studies in Mathematics, vol.~14, Walter de
  Gruyter \& Co., Berlin, 1991. \MR{MR1133391 (93e:60107)}

\bibitem{conway:90:cif}
John~B. Conway, \emph{{A course in functional analysis. 2nd ed.}}, {Graduate
  Texts in Mathematics, 96. New York etc.: Springer-Verlag. xvi, 399 p. DM
  148.00 }, 1990.

\bibitem{diestel:06:sti}
Geoff Diestel and Alexander Koldobsky, \emph{{Sobolev spaces with only trivial
  isometries.}}, Positivity \textbf{10} (2006), no.~1, 135--144.

\bibitem{doob:01:cpt}
Joseph~L. Doob, \emph{{Classical potential theory and its probabilistic
  counterpart. Reprint of the 1984 edition.}}, {Classics in Mathematics.
  Berlin: Springer. xxiii, 846 p.}, 2001.

\bibitem{federer:72:tls}
Herbert Federer and William~P. Ziemer, \emph{{The Lebesgue set of a function
  whose distribution derivatives are $p$-th power summable.}}, Math. J.,
  Indiana Univ. \textbf{22} (1972), 139--158.

\bibitem{feldman:86:obl}
W.~A. Feldman and J.~F. Porter, \emph{Operators on {B}anach lattices as
  weighted compositions}, J. London Math. Soc. (2) \textbf{33} (1986), no.~1,
  149--156. \MR{MR829395 (87j:47054)}

\bibitem{fukushima:94:df}
Masatoshi Fukushima, Y{\=o}ichi {\=O}shima, and Masayoshi Takeda,
  \emph{Dirichlet forms and symmetric {M}arkov processes}, de Gruyter Studies
  in Mathematics, vol.~19, Walter de Gruyter \& Co., Berlin, 1994.
  \MR{MR1303354 (96f:60126)}

\bibitem{goldshtejn:84:tps}
V.M. Gol'dshtejn and A.S. Romanov, \emph{{Transformations that preserve Sobolev
  spaces.}}, Sib. Math. J. \textbf{25} (1984), 382--388.

\bibitem{jost:05:pma}
J{\"u}rgen Jost, \emph{Postmodern analysis}, third ed., Universitext,
  Springer-Verlag, Berlin, 2005. \MR{MR2166001 (2006c:46001)}

\bibitem{lamperti:58:icf}
John Lamperti, \emph{On the isometries of certain function-spaces}, Pacific J.
  Math. \textbf{8} (1958), 459--466. \MR{MR0105017 (21 \#3764)}

\bibitem{lumer:63:orlicz}
Gunter Lumer, \emph{On the isometries of reflexive {O}rlicz spaces}, Ann. Inst.
  Fourier (Grenoble) \textbf{13} (1963), 99--109. \MR{MR0158259 (28 \#1485)}

\bibitem{ziemer:97:fr}
Jan Mal{\'y} and William~P. Ziemer, \emph{Fine regularity of solutions of
  elliptic partial differential equations}, Mathematical Surveys and
  Monographs, vol.~51, American Mathematical Society, Providence, RI, 1997.
  \MR{MR1461542 (98h:35080)}

\bibitem{mazya:85:tom}
V.G. Maz'ya and T.O. Shaposhnikova, \emph{{Theory of multipliers in spaces of
  differentiable functions.}}, {Monographs and Studies in Mathematics, 23.
  Pitman Advanced Publishing Program. Boston - London - Melbourne: Pitman
  Publishing Inc. XIII, 344 p. }, 1985.

\bibitem{mazya:85:ssp}
Vladimir~G. Maz{'}ya, \emph{Sobolev spaces}, Springer Series in Soviet
  Mathematics, Springer-Verlag, Berlin, 1985, Translated from the Russian by T.
  O. Shaposhnikova. \MR{87g:46056}

\bibitem{meyer-nieberg:91:bl}
Peter Meyer-Nieberg, \emph{Banach lattices}, Universitext, Springer-Verlag,
  Berlin, 1991. \MR{MR1128093 (93f:46025)}

\bibitem{meyers:75:cpp}
Norman~G. Meyers, \emph{{Continuity properties of potentials.}}, Duke Math. J.
  \textbf{42} (1975), 157--166.

\bibitem{munkres:00:top}
James~R. Munkres, \emph{{Topology. 2nd ed.}}, {Upper Saddle River, NJ: Prentice
  Hall. xvi, 537 p. }, 2000.

\bibitem{royden:88:ran}
H.L. Royden, \emph{{Real analysis. 3rd ed.}}, {New York: Macmillan Publishing
  Company; London: Collier Macmillan Publishing. xx, 444 p. }, 1988.

\bibitem{schaefer:74:blp}
Helmut~H. Schaefer, \emph{Banach lattices and positive operators},
  Springer-Verlag, New York, 1974. \MR{MR0423039 (54 \#11023)}

\bibitem{stone:37:atb}
M.~H. Stone, \emph{Applications of the theory of {B}oolean rings to general
  topology}, Trans. Amer. Math. Soc. \textbf{41} (1937), no.~3, 375--481.
  \MR{MR1501905}

\bibitem{ziemer:89:wdf}
William~P. Ziemer, \emph{Weakly differentiable functions}, Graduate Texts in
  Mathematics, vol. 120, Springer-Verlag, New York, 1989. \MR{91e:46046}

\end{thebibliography}
\bibliographystyle{amsplain}

\end{document}